\documentclass{amsart}
\usepackage{amssymb, latexsym}

\def\C{{\mathbb C}}
\def\R{{{\mathbb R}}}

\let\Re=\undefined\DeclareMathOperator*{\Re}{Re}
\let\Im=\undefined\DeclareMathOperator*{\Im}{Im}

\DeclareMathOperator*{\dist}{dist}
\DeclareMathOperator{\sech}{sech}

\DeclareMathOperator*{\Mass}{Mass}

\newcommand{\lqlr}{{L^q_t L^r_x}}
\newcommand{\ir}{{I \times \R^n}}

\newcommand{\rr}{{\R \times \R^n}}
\newcommand{\lli}{{\lqlr (\ir)}}
\newcommand{\llr}{{\lqlr (\rr)}}

\newcommand{\dnnt}{{\frac {2n}{n-2}}}

\theoremstyle{plain}
\newtheorem{theorem}{Theorem}
\newtheorem{definition}[theorem]{Definition}
\newtheorem{proposition}[theorem]{Proposition}
\newtheorem{lemma}[theorem]{Lemma}
\newtheorem{corollary}[theorem]{Corollary}

\theoremstyle{definition}
\newtheorem{remark}[theorem]{Remark}
\newtheorem*{convent}{Convention}

\numberwithin{equation}{section} \numberwithin{theorem}{section}

\begin{document}

\title{Energy-critical NLS with quadratic potentials}
\author{Rowan Killip}
\address{UCLA, Los Angeles CA 90095}
\author{Monica Visan}
\address{Institute for Advanced Study, Princeton NJ 08540}
\author{Xiaoyi Zhang}
\address{Academy of Mathematics and System Sciences, Chinese Academy of Sciences}
%\subjclass{35Q55} \keywords{nonlinear Schr\"odinger equation,
%well-posedness}

\vspace{-0.3in}
\begin{abstract}
We consider the defocusing $\dot H^1$-critical nonlinear Schr\"odinger equation in all dimensions
($n\geq 3$) with a quadratic potential $V(x)=\pm \tfrac12 |x|^2$.  We show global well-posedness
for radial initial data obeying $\nabla u_0(x),\  xu_0(x) \in L^2$.  In view of the potential
$V$, this is the natural energy space.  In the repulsive case, we also prove scattering.

We follow the approach pioneered by Bourgain and Tao in the case of no potential; indeed,
we include a proof of their results that incorporates a couple of simplifications
discovered while treating the problem with quadratic potential.
\end{abstract}

\maketitle

\tableofcontents

\section{Introduction}

We study the Cauchy problem for the defocusing energy-critical
nonlinear Schr\"odinger equation with potential,
\begin{equation}\label{equation}
\begin{cases}
\ iu_t = -\frac12 \Delta u +  V u+|u|^{\frac{4}{n-2}}u, \\
\ u(0)=u_0
\end{cases}
\end{equation}
in the cases $V(x)=\pm \tfrac12 |x|^2$.  Here $u(t,x)$ is a complex-valued function of spacetime $\R\times \R^n$ and $n\geq3$.
This is a Hamiltonian PDE with energy function
\begin{equation}\label{energy}
E(u(t)) := \int_{\R^n} \tfrac 12|\nabla u(t,x)|^2 + V(x)|u(t,x)|^2 + (1-\tfrac2n)|u(t,x)|^{\dnnt} \, dx.
\end{equation}
In particular, the energy is conserved as explained below.

It is natural to consider \eqref{equation} for arbitrary finite-energy initial data $u_0$.
In the case of the confining potential, $V(x)=\tfrac12 |x|^2$, all terms in \eqref{energy}
have the same sign and so the energy class is
$$
\Sigma := \{ f : \|f\|_{\Sigma}^2:= \|\nabla f\|_2^2 + \|x f\|_2^2 < \infty \}.
$$
(Note that the $\Sigma$-norm controls the term $\int_{\R^n} |u(t,x)|^{\dnnt}\,dx$ in \eqref{energy} by Sobolev embedding.)

In the repulsive case, $V(x)=- \tfrac12 |x|^2$, things are more subtle.  To us, the natural choice is to require
that both the kinetic energy
\begin{equation}\label{K energy}
\tfrac 12 \int_{\R^n} |\nabla u(t,x)|^2 \, dx
\end{equation}
and the potential energy
\begin{equation}\label{P energy -}
\int_{\R^n} -\tfrac12 |x|^2 |u(t,x)|^2 + (1-\tfrac2n)|u(t,x)|^{\dnnt} \, dx
\end{equation}
are finite.  By Sobolev embedding, $\Sigma$ is exactly the class of initial data satisfying these two
requirements.

The name energy-critical relates to the power $\frac{4}{n-2}$ appearing in \eqref{equation}.  If we were to discard
the potential for a moment, the scaling
\begin{equation}\label{scaling}
u(t,x) \mapsto u_\lambda (t,x) := \lambda^{\frac{2-n}{2}} u\bigl(\lambda^{-2}t,\lambda^{-1}x\bigr)
\end{equation}
maps a solution to \eqref{equation} to another solution to \eqref{equation}.  Moreover, this scaling also
leaves the energy invariant; hence the name energy-critical.  Put more simply, the potential and kinetic energies
scale in the same way and so are comparable at all length scales.  In subcritical problems, the kinetic energy
dominates the potential energy at small length scales; in supercritical problems the roles are reversed and proving
well-posedness becomes very difficult indeed.

With the introduction of a quadratic potential, the equation no longer has a scaling symmetry; however, it
is still natural to consider $\frac{4}{n-2}$ as the critical power.  In particular, the dichotomy between
blowup and well-posedness occurs at very small length scales, where one may approximate the potential by
a constant.

Our goal here is to construct global \emph{strong solutions} to \eqref{equation}, that is, $u\in C^0(\R;\Sigma)$
obeying the Duhamel formula
\begin{equation}\label{duhamel}
u(t)=U(t-t_0)u(t_0) - i\int_{t_0}^t U(t-s) \bigl(i\partial_s+\tfrac 12 \Delta -V\bigr) u(s)\, ds.
\end{equation}
Here $U(t)$ denotes the linear propagator associated to \eqref{equation}, $U(t)=e^{it(\frac12\Delta-V)}$.
As explained in \cite[\S 3.3]{cazenave:book}, energy is conserved for such solutions.  In this paper, we
achieve this goal for radial initial data:

\begin{theorem}\label{confining gwp}
Let $V(x)=\tfrac 12 |x|^2$ and $u_0\in \Sigma$ be radial. Then, there exists a unique global strong solution $u$ to \eqref{equation}
and on each compact time interval $[-T,T]$,
\begin{equation}\label{TE:conf}
\|H_0^{\frac 12} u\|_{ S^0([-T,T])}\le C(T,\|u_0\|_\Sigma).
\end{equation}
\end{theorem}

\begin{theorem}\label{repulsive gwp}
Let $V(x)=-\tfrac 12|x|^2$ and $u_0\in \Sigma$ be radial. Then, there exists a unique global strong solution $u$ to \eqref{equation}.
Moreover,
\begin{equation}\label{TE:repul}
\|H(-t)^{\frac 12} u\|_{ S^0(\R)}\le C(\|u_0\|_\Sigma)
\end{equation}
and there exist unique radial functions $u_{\pm}\in \Sigma$ such that
$$
\|U(-t)u(t)-u_{\pm}\|_{\Sigma}\to 0, \quad  \text{as } t\to \pm\infty.
$$
\end{theorem}

The $S^0$ norm is the supremum over a natural class of spacetime norms and $H_0=-\frac12\Delta + \tfrac12 |x|^2$
is the harmonic oscillator.  The operator $H(-t)$ is an analogue of $H_0$ adapted to the repulsive case.  For
a full description of our notations, see Section~\ref{not}.  Let us remark here however that \eqref{TE:conf}
and \eqref{TE:repul} are strong enough to bound the $L_{t,x}^{{2(n+2)}/{(n-2)}}$ norm of the solution;
see \eqref{emb} in the confining case and Lemma~\ref{stric bdd} in the repulsive case.

We will now give a brief survey of related works, emphasizing the origins of some of the techniques we use.

The paper \cite{YGOh} treats the nonlinear Schr\"odinger equation with general non-negative potential and $L^2$-subcritical power nonlinearity.
Global well-posedness is proved provided the potential obeys $\partial^\alpha V \in L^\infty$ for $|\alpha|\geq 2$.

In the mathematical community, the most active proponent of NLS with quadratic potential has been Carles.  In particular, he
has shown global well-posedness for energy-subcritical power nonlinearity; see \cite{carles1} for the confining case,
\cite{carles2} for the repulsive case, and \cite{carles3} in the sign-indefinite case.  We will use one of the techniques he
introduced in the repulsive case.  By splitting the energy in a carefully chosen way, he was able to obtain \emph{a priori}
exponential decay (cf. \eqref{pot decay}).  This is vital in the proof of scattering.

In the physical literature, \eqref{equation} with confining potential has been used to describe Bose--Einstein condensates
in a trap.  See for example, the discussion of the Gross--Pitaevskii equation in the lecture course \cite{CoT}.

We should also mention that there has been a wealth of work on the focusing problem, that is, when the nonlinearity has the
opposite sign.  The main interest in that case has been the blowup of solutions.  See the overview \cite{carlesJEDP} for an
introduction to this field.

Thus far, we have described results for energy-subcritical nonlinear Schr\"odinger equation with potential.  The second thread
leading to this paper has been the development of an energy-critical theory for NLS without potential.  Indeed, the passage from
subcritical to critical nonlinearities had been a major hurdle in the theory until breached by Bourgain, \cite{borg:scatter}.

Bourgain proved global well-posedness and scattering for the three dimensional equation with radial data.
A different approach to global regularity was subsequently developed by Grillakis, \cite{grillakis:scatter},
again for radially symmetric initial data.

In \cite{zhang}, Zhang adapted Bougain's argument to the case of a repulsive potential, proving global well-posedness and
scattering for \eqref{equation} in three dimensions.

Returning to our discussion of the case $V\equiv 0$, Bourgain's argument was simplified and extended to all dimensions ($n\geq 3$)
by Tao, \cite{tao: gwp radial}.  We will prove Theorems~\ref{confining gwp} and~\ref{repulsive gwp} by mimicking the argument given
there. In particular, our proof in the three-dimensional case is significantly simpler than that in \cite{zhang}.

In adapting the argument from \cite{tao: gwp radial}, we uncovered some simplifications.  As a consequence, we present a still simpler
proof of the main results of \cite{borg:scatter,tao: gwp radial}:

\begin{theorem}\label{free gwp}
Let $V\equiv 0$ and $u_0\in \dot H^1(\R^n)$ be radial. Then, there exists a unique global strong solution $u$ to \eqref{equation}
and
$$
\|\nabla u\|_{ S^0(\R)}\le C(E(u_0)).
$$
Moreover, there exist unique radial functions $u_{\pm}\in \dot H^1(\R^n)$ such that
$$
\|e^{-it\Delta/2} u(t)-u_{\pm}\|_{\dot H^1(\R^n)}\to 0, \quad  \text{as } t\to \pm\infty.
$$
\end{theorem}

Our main simplification relative to \cite{tao: gwp radial} is the improved perturbation theory discussed in Section~\ref{S:pert};
one should compare our Lemma~\ref{pert lemma} with Tao's Appendix.  This also leads to a cleaner treatment of the concentration result,
Proposition~\ref{bubble}.

The roots of our perturbation lemma lie in \cite{TV} (see \cite{nak} for a similar approach to the
Klein--Gordon equation).  The idea is to work in spaces with a fractional number of derivatives but critical scaling.  Tao
works in a space with no derivatives (and critical scaling); this is well adapted for proving the concentration result, but
results in gymnastics for the perturbation lemma in high dimensions, $n\geq 6$.  If one were to work in a space with a full
derivative, the perturbation theory would become trivial, while making the concentration result much more difficult.
The middle road, a fractional number of derivatives, simplifies the proofs of both the concentration and perturbation results.

Working in an exotic space requires exotic estimates.  These are of two types.  Firstly, we need the exotic Strichartz estimates
of Foschi, \cite{foschi}, which hold under very mild hypotheses.

Secondly, we need a fractional chain rule, such as that proved by Christ and Weinstein, \cite{christw}.  However, because of
the need to commute differentiation with the linear propagator, the `derivatives' that appear in the spaces we have been discussing
are not fractional powers of the Laplacian, but rather fractional powers of different (time-dependent) second order differential operators.
The specific fractional chain rule we need is given in Proposition~\ref{frac chain}.

One of the main tools in the treatment of critical problems is the Morawetz inequality. It is a type of monotonicity formula
related to the conservation of momentum, which corresponds to $\dot H^{\frac 12}$-scaling.  As noticed already by Bourgain,
one can recover energy-critical bounds by restricting it to finite volume, thus suppressing the low-frequency problem.

The particular Morawetz inequality we use is derived in Proposition~\ref{morawetz}.  As it prevents concentration at the spatial origin,
it is particularly well adapted to the treatment of spherically symmetric solutions.  These can only concentrate significantly
near the origin; see Corollary~\ref{O bubble}.

An interaction form of Morawetz inequality was introduced in \cite{ckstt:?}; it has permitted the treatment of energy-critical
problems without the assumption of radial symmetry, \cite{ckstt:gwp}.  This was a major breakthrough and has now been
extended to all dimensions, \cite{RV, MV}.  Adapting the \emph{a priori} interaction inequality (which scales like
$\dot H^{\frac14}$) to the energy-critical setting is a major undertaking and requires the introduction of frequency cutoffs.
Implementing such a strategy in the presence of a potential would constitute an ambitious project.

The paper is composed as follows.  In Section~2 we introduce notation and the basic estimates we will use.  In particular,
the fractional chain rule is developed there.  The Morawetz inequality is derived in Section~3 along with bounds on
the transportation of mass.  In Section~4 we discuss the local theory for \eqref{equation}; this is mostly a repetition of
well-known techniques.  Section~5 treats perturbation theory; this is the main simplification compared to \cite{tao: gwp radial}.
With the exception of scattering for the repulsive potential, the theorems stated above are proved in Section~6.  In Section~7
we show existence of wave operators and asymptotic completeness in the repulsive case.

\textbf{Acknowledgements}: The first author was supported in part by NSF grant DMS--0401277 and a Sloan Foundation Fellowship.
The work of the second author was supported in part by the NSF under grant DMS--0111298.
The third author was supported by the NSF grant No. 10601060 (China).  She would also like to thank the hospitality of UCLA,
where much of this work was completed.

We are grateful to the referee for several expository remarks and to Dong Li for pointing out an error in a previous proof of Lemma~\ref{L:boho}.

Any opinions, findings and conclusions or recommendations expressed in this material
are those of the authors and do not reflect the views of the National Science Foundation.

%%%%%%%%%%%%%%%%%%%%%%%%%%%%%%%%%%%%%%%%%%%%%%%%%%%%%%%%%%%%%%%%%%%%%%%%%%%%%%%%%%%%%%%%%%%
%
%
%                                   Section
%
%
%%%%%%%%%%%%%%%%%%%%%%%%%%%%%%%%%%%%%%%%%%%%%%%%%%%%%%%%%%%%%%%%%%%%%%%%%%%%%%%%%%%%%%%%%%%

\section{Notations and basic estimates}\label{not}
We will often use the notations $X\lesssim Y$ and $X=O(Y)$ to mean that there exists some constant so that $X\le CY$. Similarly, we will use
$X\simeq Y$ if $X\lesssim Y\lesssim X$. The derivative operator $\nabla$ refers to the space variable only. We will occasionally use subscripts
to denote spatial derivatives and will use the summation convention over repeated indices.

We use $U(t)$ to denote the linear Schr\"odinger propagator associated with \eqref{equation}, that is, $U(t)=e^{it(\frac 12 \Delta-V)}$.
When $V\equiv 0$, this is the free Schr\"odinger propagator whose explicit expression can be found via the Fourier transform:
$$
[U(t)f](x) = \frac 1{(2\pi i t)^{\frac n2}}\int_{\R^n} e^{\frac{i(x-y)^2}{2t}} f(y)\, dy.
$$
When $V$ is the harmonic potential, $U(t)$ is given by Mehler's formula, \cite{Folland}, as follows:
\begin{align*}
[U(t)f](x)&=\frac 1{(2\pi i\sinh t)^{\frac n2}}\int_{\R^n} e^{\frac i{\sinh t}(\frac {x^2+y^2}2\cosh t-xy)}f(y)\, dy,
        \quad \text{for } V(x)=-\tfrac 12 |x|^2,\\
[U(t)f](x)&=\frac 1{(2\pi i\sin t)^{\frac n2}}\int_{\R^n} e^{\frac i{\sin t}(\frac {x^2+y^2}2\cos t-xy)}f(y)\, dy,
        \quad \text{for } V(x)=\tfrac 12 |x|^2.
\end{align*}
In particular, in the free and the repulsive potential cases the propagator obeys the dispersive estimate
\begin{align}\label{disp est}
\|U(t)f\|_{L_x^{\infty}}\lesssim |t|^{-\frac n2}\|f\|_{L_x^1}, \quad \text{for all } t\neq 0.
\end{align}
A similar estimate holds in the case of the confining potential; more precisely,
\begin{align}\label{conf disp est}
\|U(t)f\|_{L_x^{\infty}}\lesssim |\sin t|^{-\frac n2}\|f\|_{L_x^1}, \quad \text{for all } t\neq 0,
\end{align}
which yields \eqref{disp est} when restricted to $t\in [-\pi/2, \pi/2]$.

We use $L_x^r(\R^n)$ to denote the Banach space of functions $f:\R^n\to \C$ whose norm
$$
\|f\|_r:=\Bigl(\int_{\R^n} |f(x)|^r dx\Bigr)^{\frac{1}{r}}
$$
is finite, with the usual modifications when $r=\infty$.  We use $\lqlr$ to denote the spacetime norm
$$
\|u\|_{\llr} :=\Bigl(\int_{\R}\Bigl(\int_{\R^n} |u(t,x)|^r dx \Bigr)^{q/r} dt\Bigr)^{1/q},
$$
with the usual modifications when $q$ or $r$ is infinity, or when the domain $\R \times \R^n$ is
replaced by some smaller spacetime region.  When $q=r$ we abbreviate $\lqlr$ by $L^q_{t,x}$.

We say that a pair of exponents $(q,r)$ is \emph{Schr\"odinger admissible} if $\tfrac{2}{q} + \tfrac{n}{r} = \frac{n}{2}$ with $2 \leq q,r \leq \infty$
and $(q,r,n)\neq (2,\infty,2)$.  If $I \times \R^n$ is a spacetime slab, we define the \emph{Strichartz norm} $S^0(I)$ by
\begin{equation}\label{s0}
\|u\|_{S^0(I)} := \sup \|u \|_{\lli}
\end{equation}
where the $\sup$ is taken over all Schr\"odinger-admissible pairs $(q,r)$.
We call this a Strichartz norm due to its appearance in the inequalities of the same name:

\begin{lemma}\label{strichartz estimate}
Let $u:I\times \R^n$ be an $S^0$ solution to the Schr\"odinger equation
$$
iu_t+\tfrac12 {\Delta} u-V  u=\sum_{i=1}^M F_m,
$$
for some functions $F_1,\ldots, F_m$.  Then,
$$
\|u\|_{S^0(I)}
\lesssim \|u(t_0)\|_2+\sum_{m=1}^M\|F_m\|_{L_t^{q_m'}L_x^{r_m'}(I\times\R^n)}
$$
for any $t_0\in I$ and admissible pairs $(q_m,r_m)$. As usual, $p'$ denotes the dual exponent to $p$, $\frac 1p+\frac 1{p'}=1$.
The implicit constant is absolute except in the confining case where is depends on $|I|$.
\end{lemma}

\begin{proof}
Strichartz first proved an inequality of this type for $V\equiv 0$, see \cite{strichartz}.  The result for the full range of exponents was
proved by Keel and Tao, \cite{tao:keel}.  Indeed, they show that these estimates follow from just two hypotheses: conservation
of the $L^2_x$-norm and the dispersive estimate \eqref{disp est}, both of which hold for the three propagators we consider here.
\end{proof}

From Sobolev embedding, we have

\begin{lemma}\label{lemma str-norms}
For any function $u$ on $\ir$ we have
\begin{align}\label{str-norms}
\|\nabla u\|_{L^\infty_t L^2_x}
& + \|\nabla u \|_{L^{\frac{2(n+2)}{n-2}}_t L^{\frac{2n(n+2)}{n^2+4}}_x} + \|\nabla u\|_{L^{\frac{2(n+2)}{n}}_{t,x}}
        + \|\nabla u\|_{L^2_t L^{\frac{2n}{n-2}}_x} \notag\\
& + \|u\|_{L^\infty_t L^{\frac{2n}{n-2}}_x} + \|u\|_{L^{\frac{2(n+2)}{n-2}}_{t,x}}
        + \|u\|_{L^{\frac{2(n+2)}{n}}_t L^{\frac{2n(n+2)}{n^2-2n-4}}_x} \lesssim \|\nabla u\|_{S^0},
\end{align}
where all spacetime norms are on $\ir$.
\end{lemma}

To keep the formulae readable, we introduce several abbreviated notations.  For a time interval $I$ we define
$$
\|u \|_{Z(I)}:= \|u\|_{L^{\frac{2(n+2)}{n-2}}_{t,x}(I\times\R^n)}.
$$
As described in the introduction, one of the main simplifications of this paper relative to \cite{tao: gwp radial} is the decision
to work with spacetime norms with a fractional number of derivatives but invariant under the scaling \eqref{scaling}.
In the free case, the meaning of derivatives is clear; in the confining case, $V=\frac12|x|^2$, it will be based off the Hermite operator
$$
H_0 = -\tfrac 12 \Delta + \tfrac 12 |x|^2.
$$
In the repulsive case, the choice is more subtle: $H(t) = U(-t) H_0 U(t)$.  To unify the notation we make the following
general definition:
\begin{equation*}
H(t) = \begin{cases}
\ -\Delta           &\quad \text{for }  V\equiv 0 \\
\ H_0               &\quad \text{for } V(x)=\tfrac12 |x|^2\\
\ U(-t) H_0 U(t)    &\quad \text{for } V(x)=-\tfrac12 |x|^2.
\end{cases}
\end{equation*}
With this in place, we define our main spacetime norm, $W$,
\begin{equation}\label{W defn}
\|u\|_{W(I)}:=  \| H(-t)^{\frac 12 -\frac 1n} u \|_{L_t^{\frac{2(n+2)}{n-2}}L_x^{\frac{2n^2(n+2)}{n^3-8}}(I\times\R^n)}.
\end{equation}
Note that here, as elsewhere in the paper, functions of $H(-t)$ are defined via the $L^2$ functional calculus.
By Lemma~\ref{lemma str-norms} and Sobolev embedding, in the free case we have
$$
\|u\|_{Z(I)} \lesssim \|u\|_{W(I)} \lesssim \| H(-t)^{\frac12} u\|_{ S^0(I)}.
$$
In fact, this also holds in the presence of a quadratic potential; see Lemma~\ref{sob em}.

In order to estimate the $W$-norm of the solution, we will need to estimate the nonlinearity
in a corresponding space, $N$.   We choose to define this as follows:
\begin{equation}\label{N defn}
\| F \|_{N(I)}:= \| H(-t)^{\frac 12 -\frac 1n}  F\|_{L_t^{2}L_x^{\frac{2n^2}{n^2+2n-4}}(I\times\R^n)}.
\end{equation}
Using the Strichartz inequalities of Foschi, \cite{foschi}, we know that for $n\geq 5$,
\begin{equation}\label{foschi est}
\Bigl\| \int_{t_0}^t U(t-s) F(s) \,ds \Bigr\|_{W(I)} \lesssim \| F \|_{N(I)}
\end{equation}
for all $t_0\in I$.  Because of the poor dispersive estimate for the confining potential, the implicit constant
in \eqref{foschi est} depends on $|I|$ in this case.

\subsection{Heisenberg's equations}  In this subsection, we discuss the solutions of Heisenberg's equations
for the evolution of the momentum and position operators:
\begin{align*}
\tfrac{d\ }{dt} P(t) &= i[-\tfrac12\Delta + V , P(t)] \qquad P(0) = i\nabla \\
\tfrac{d\ }{dt} X(t) &= i[-\tfrac12\Delta + V , X(t) ] \qquad X(0) = x.
\end{align*}
When $V(x)=-\tfrac 12 |x|^2$, the solutions are
\begin{align}\label{p of t -}
P(t) &= U(-t)i\nabla U(t)= i\nabla \cosh t + x\sinh t \\
    &= i\cosh(t) \, e^{-\frac{i}2 |x|^2\tanh t}\nabla e^{\frac{i}2 |x|^2\tanh t} \notag \\
\label{x of t -}
X(t) &= U(-t)x U(t)= i\nabla \sinh t + x\cosh t \\
     &= i\sinh(t) \, e^{-\frac{i}2 |x|^2\coth t}\nabla e^{\frac{i}2 |x|^2\coth t}, \notag
\end{align}
while for $V(x)=\tfrac 12 |x|^2$ one obtains
\begin{align}\label{p of t +}
P(t) &= U(-t)\,i\nabla\, U(t)= i\nabla \cos t - x\sin t= i\cos(t) e^{-\frac{i}2 |x|^2\tan t}\nabla e^{\frac{i}2 |x|^2\tan t}\\
\label{x of t +}
X(t) &= U(-t)\,x\,U(t)= i\nabla \sin t + x\cos t = i\sin(t) e^{\frac{i}2 |x|^2\cot t}\nabla e^{-\frac{i}2 |x|^2\cot t}.
\end{align}

In the repulsive case, $V=-\tfrac12 |x|^2$, we may invert the equations above to obtain
\begin{align}\label{x and deriv}
i\nabla=\cosh(t) P(t) - \sinh(t) X(t),\ x=\cosh(t) X(t)-\sinh(t) P(t).
\end{align}
A similar result holds in the confining case,
$$
i\nabla =\cos(t) P(t) + \sin(t) X(t),\ x=\cos(t) X(t) - \sin(t) P(t);
$$
however, we will not need these formulae.

As noted in \cite{carles1}, $P(t)$ and $X(t)$ behave almost like derivatives in the sense described in the
next two lemmas.

\begin{lemma}\label{lemma jt ht}
For $V(x)=\pm \frac12  |x|^2$ and any nonlinearity $F(z)=|z|^p z$, $p>0$, we have
\begin{align*}
P(t)F(u)&=\partial_z F(u) P(t)u  - \partial_{\bar z} F(u)\overline{P(t)u}\\
X(t)F(u)&=\partial_z F(u) X(t)u  - \partial_{\bar z} F(u)\overline{X(t)u},
\end{align*}
where  $P(t)$ and $X(t)$ are as defined above.
\end{lemma}

\begin{proof}
This result is not peculiar to $P(t)$ and $X(t)$; it holds for any operator of the form
$a(x)i\nabla + b(x)$ with real-valued functions $a$ and $b$.  On the one hand,
%; it holds for any linear combination of the operators $i\nabla$ and $x$.  Specifically,
$$
i\nabla F(u)= \partial_z F(u) i\nabla u  + \partial_{\bar z} F(u) i \nabla \bar u
$$
by the chain rule.  For multiplication operators, one should simply note that
$$
F(u) = \partial_z F(u) u  - \partial_{\bar z} F(u) \bar{u}
$$
due to the specific form of the nonlinearity.
\end{proof}

\begin{lemma}[Sobolev embedding]\label{lemma pt sobolev}
For exponents $1<p<p^*<\infty$ obeying $\frac 1{p^*}=\frac 1{p}-\frac 1n,$  we have
\begin{align*}
\| f \|_{p^*} \lesssim  \frac{1}{\cosh(t)} \| P(t) f \|_{p}
\end{align*}
when $V(x)=- \frac12 |x|^2$ and
\begin{align*}
\| f \|_{p^*} \lesssim  \frac{1}{|\cos(t)|} \| P(t) f \|_{p}
\end{align*}
when $V(x)= \frac12 |x|^2$.
\end{lemma}

\begin{proof}
This follows from the usual Sobolev embedding together with the right-most formulae in
\eqref{p of t -} and \eqref{p of t +}.
\end{proof}

\subsection{Energy control}\label{energy control ss}
When $V\equiv 0$, conservation of the energy, \eqref{energy}, directly gives
\begin{align}\label{H(-t) control free}
\|\nabla u(t)\|_2\lesssim E^{\frac 12}.
\end{align}
Similarly, when $V(x)=\frac12 |x|^2$, we have
\begin{align}\label{H(-t) control conf}
\| u(t)\|_\Sigma \lesssim E^{\frac12}.
\end{align}
However, when $V(x)=-\frac12|x|^2$ energy conservation does not yield control on $\|\nabla u(t)\|_2$ or $\|xu(t)\|_2$.
Therefore, we adopt the approach of \cite{carles1}, namely to split the energy into two parts,
\begin{align*}
\mathcal E_1(t)=\tfrac12\|P(-t)u(t)\|_2^2+(1-\tfrac2n)\cosh^2(t) \|u(t)\|_{\dnnt}^{\dnnt},\\
\mathcal E_2(t)=\tfrac12\|X(-t)u(t)\|_2^2+(1-\tfrac2n)\sinh^2(t) \|u(t)\|_{\dnnt}^{\dnnt}.
\end{align*}
Notice that each is positive and $E(u(t)) = \mathcal E_1(t) - \mathcal E_2(t)$.
It is not hard to verify the following analogue of the variance identity
\begin{align}\label{deriv eps}
\frac{d\mathcal E_1(t)}{dt} = \frac{d\mathcal E_2(t)}{dt} =  - \frac2n \sinh(2t) \|u(t)\|_{\dnnt}^{\dnnt};
\end{align}
for details, see \cite{carles1}.  As in that paper, this leads to

\begin{lemma}[Energy control in the repulsive case]\label{lemma energy control}
The potential energy decays exponentially in time,
\begin{align}\label{pot decay}
 \|u(t)\|_{\dnnt}^{\dnnt} \lesssim \mathcal E_1(0) \cosh^{-2}(t),
\end{align}
while the other parts of the energy do not grow too rapidly,
\begin{align}
\|H(-t)^\frac12 u(t)\|_2 &\lesssim C(\|u_0\|_\Sigma),           \label{H(-t) control repuls} \\
\|u(t)\|_\Sigma &\lesssim \cosh(t) C(\|u_0\|_\Sigma).   \label{sigma control}
\end{align}
\end{lemma}

\begin{proof}
From \eqref{deriv eps} we see that $\mathcal E_1$ achieves its maximum at $t=0$, which immediately yields \eqref{pot decay}.
Similarly,
$$
2 \|H(-t)^\frac12 u(t)\|_2^2 = \|P(-t)u(t)\|_2^2 + \|X(-t)u(t)\|_2^2 \lesssim \mathcal E_1(0) + \mathcal E_2(0) \lesssim C(\|u_0\|_\Sigma),
$$
where the second inequality follows from Sobolev embedding.  From this and \eqref{x and deriv} we obtain \eqref{sigma control}.
\end{proof}

\subsection{Fractional chain rule for the harmonic oscillator}  The goal of this subsection is to prove a form
of fractional chain rule.  The main application of this will be to bound the solution in the space $W$. Note that
in the repulsive case this forces us to consider
\begin{align*}
H(-t) &= U(t) H_0 U(-t) = \tfrac12 P(-t)^2 + \tfrac12 X(-t)^2 \\
&= \cosh(2t) H_0 - \tfrac i2 \sinh(2t) (x\cdot\nabla + \nabla\cdot x).
\end{align*}

\begin{lemma}\label{L:H one way}
Let $H=-\tfrac{a}2\Delta + \tfrac{ib}2(x\cdot\nabla + \nabla\cdot x) + \tfrac{c}2 |x|^2$
with $a,c>0$ and normalized so that $ac-b^2=1$.  Then for all $1<p<\infty$,
\begin{align}\label{H est}
\bigl\| -\Delta f \bigr\|_{p} + \bigl\| (x\cdot\nabla + \nabla\cdot x) f \bigr\|_{p} + \bigl\| |x|^2 f \bigr\|_{p}
\simeq \bigl\| H f \bigr\|_{p},
\end{align}
where the implicit constant depends continuously on $a$, $b$, and $c$.
\end{lemma}

\begin{proof}
The $\gtrsim$ portion of \eqref{H est} follows from the triangle inequality, so we need only consider the
other inequality.

We begin with the special case $H_0=-\tfrac12\Delta + \tfrac12 |x|^2$.  By Mehler's formula, the operator $e^{-tH_0}$
has Weyl symbol
$$
 [\sech(\tfrac{t}2)]^{n} \exp\bigl\{ - \tanh(\tfrac{t}2)[4\pi^2\xi^2+|x|^2] \bigr\};
$$
see \cite[Ch. 5]{Folland}.  The identity $H_0^{-1}=\int_0^\infty e^{-tH_0}\,dt$ allows us to deduce the
symbol for the resolvent:
$$
\sigma(x,\xi) =  2\int_0^1 \exp\bigl\{ - [4\pi^2\xi^2+|x|^2] u \bigr\} (1-u^2)^{\frac{n}2-1}\,du
$$
where $u=\tan(t/2)$.  This makes it easy to check that the symbols for $|x|^2 H_0^{-1}$, $-\Delta H_0^{-1}$,
and $(x\cdot\nabla + \nabla\cdot x) H_0^{-1}$ belong to the symbol class $S^0$.
Therefore, we may conclude that \eqref{H est} holds for $H_0$; see \cite[\S VI.5]{stein:large}.

Now we turn to general $H$.  We will assume $a=1$, since the dependence on $a$ can be recovered by
scaling.  This permits us to use the operator identity
$$
e^{ib|x|^2/2} (-\Delta + |x|^2) e^{-ib|x|^2/2} = -\Delta + i b(\nabla\cdot x + x \cdot \nabla) + (1+b^2) |x|^2,
$$
which implies $H^{-1}=e^{ib|x|^2/2} H_0^{-1} e^{-ib|x|^2/2}$.  In this way, we may estimate
\begin{align*}
\bigl\| -\Delta H^{-1} f \bigr\|_{p}
&\leq \bigl\| -\Delta H_0^{-1} e^{-ib|x|^2/2}  f \bigr\|_{p} + b^2 \bigl\| |x|^2 H_0^{-1} e^{-ib|x|^2/2}  f  \bigr\|_{p} \\
&\qquad {} + |b|\,\bigl\| (\nabla\cdot x + x \cdot \nabla) H_0^{-1} e^{-ib|x|^2/2}  f  \bigr\|_{p} \\
&\lesssim (1+b^2)\bigl\| f \bigr\|_{p},
\end{align*}
using $\|e^{\pm ib|x|^2/2}  g\|_p = \|g\|_p$ and the $H_0$ version of \eqref{H est} in the process.
Similarly, $\bigl\||x|^2 H^{-1} f \bigr\|_p = \bigl\| |x|^2 H_0^{-1} e^{-ib|x|^2/2} f \bigr\|_p \lesssim \bigl\|f\bigr\|_p$.
Lastly,
\begin{align*}
\bigl\| (\nabla\cdot x + x \cdot \nabla) H^{-1} f \bigr\|_p
&\leq \bigl\| (\nabla\cdot x + x \cdot \nabla) H_0^{-1} e^{-ib|x|^2/2}  f \bigr\|_p \\
&\qquad + 2|b|\,\bigl\| |x|^2 H_0^{-1} e^{-ib|x|^2/2}  f  \bigr\|_p  \\
&\lesssim \sqrt{1+b^2}\,\bigl\| f \bigr\|_p,
\end{align*}
which completes the proof.
\end{proof}

\begin{lemma}\label{L:same norm}
Let $H=-\tfrac{a}2\Delta + \tfrac{ib}2(x\cdot\nabla + \nabla\cdot x) + \tfrac{c}2 |x|^2$
with $a,c>0$ and normalized so that $ac-b^2=1$.  Then
\begin{align}\label{E:same norm}
 \|H^{\gamma} f \|_{p} \simeq  \| (-\Delta)^{\gamma} f\|_{p} + \|\, |x|^{2\gamma}f\|_{p}
\end{align}
for each $\gamma\in[0,1]$ and $1< p <\infty$.  The implicit constants depend continuously on $a$, $b$, and $c$.
\end{lemma}

\begin{proof}
Naturally, we will prove this result using analytic interpolation (see \cite[\S V.4]{SteinWeiss}).
To pursue this approach, one needs a bound such as
$$
\bigl\| (-\Delta)^{is} f\bigr\|_{p} + \bigl\| |x|^{2is} f \bigr\|_{p}
    + \bigl\| H^{is} f \bigr\|_{p} \lesssim \bigl\| f \bigr\|_{p}
\quad\text{for all $s\in\R$.}
$$
For $|x|^{2is}$ this is trivial, while for $(-\Delta)^{is}$ one may use the Mikhlin multiplier theorem.
The last term is less traditional. Fortunately, Hebisch, \cite{Hebisch}, has proved a general multiplier theorem
for Schr\"odinger operators. In particular, $H_0^{is}$ is uniformly bounded on $L^p$.  This then implies the
same for $H^{is}$ via conjugation with $e^{ib|x|^2/2}$, just as in the preceding proof.

By Lemma~\ref{L:H one way} we know that
\begin{align*}
  \| (-\Delta) H^{-1} f\|_{p} + \|\,|x|^{2} H^{-1} f\|_{p} \lesssim \| f \|_{p},
\end{align*}
which corresponds to $\gamma=1$; the $\gamma=0$ case is trivial.  By employing
analytic interpolation with the operator-valued functions $z\mapsto \Delta^{z}H^{-z}$ and $z\mapsto |x|^{2z}H^{-z}$, one may deduce
\begin{align}\label{E:sn 1}
   \| (-\Delta)^\gamma f\|_{p} + \|\,|x|^{2\gamma} f\|_{p} \lesssim \|H^\gamma  f \|_{p},
\end{align}
for $\gamma\in[0,1]$, which gives us one inequality in \eqref{E:same norm}.  By the same argument using
the function $z\mapsto H^{z}H_0^{-z}$,
\begin{align}\label{E:sn 2}
  \| H^{\gamma} f\|_{p} \lesssim \|H_0^\gamma  f \|_{p},
\end{align}
which we will now use to obtain the other direction in \eqref{E:same norm}. We argue as follows
\begin{align*}
\|H^\gamma  f \|_{p} &\lesssim  \| H_0^\gamma f \|_{p} \\
&\lesssim \|H_0^{\gamma-1}(-\Delta)^{1-\gamma}\|_{p\to p} \cdot \|(-\Delta)^{\gamma} f \|_{p} + \\
    &\hspace*{20mm} {} + \|H_0^{\gamma-1}|x|^{2-2\gamma}\|_{p\to p} \cdot \|\,|x|^{2\gamma} f \|_{p} \\
&\lesssim \|(-\Delta)^{\gamma} f \|_{p} + \|\,|x|^{2\gamma} f \|_{p}.
\end{align*}
The second step is essentially the triangle inequality, while the last step consists of using the dual of \eqref{E:sn 1}
with $\gamma\mapsto1-\gamma$ and $p\mapsto p'$.
\end{proof}

Combining this lemma with the usual Sobolev embedding, we obtain an analogue for the harmonic oscillator:

\begin{lemma}[Sobolev embedding]\label{sob em} Let $H=-\tfrac{a}2\Delta + \tfrac{ib}2(x\cdot\nabla + \nabla\cdot x)
+ \tfrac{c}2 |x|^2$ with $a,c>0$ and normalized so that $ac-b^2=1$. Given $\gamma\in[0,1]$ and $p\in(1,\infty)$,
$$
\bigl\| f \bigr\|_q \lesssim \bigl\| H^\gamma f \bigr\|_p
$$
provided $\frac1q = \frac1p - \frac{2\gamma}n$.  In particular,
\begin{equation}\label{emb}
\|u\|_{Z(I)} \lesssim \|u\|_{W(I)} \lesssim \| H(-t)^{\frac12} u\|_{ S^0(I)}.
\end{equation}
\end{lemma}

\begin{remark}
Tracing the dependence on $a$, $b$, and $c$ of the implicit constant in the Sobolev embedding inequalities above,
in the repulsive case we find
\begin{equation}\label{remark2.9}
\|u\|_{Z(I)} \lesssim \cosh^{1-\frac 2n}(T) \|u\|_{W(I)} \lesssim \cosh(T)\| H(-t)^{\frac12} u\|_{ S^0(I)},
\end{equation}
where $T$ is such that $I\subseteq [-T,T]$.
\end{remark}

We are now ready to prove the main result of this subsection.  It is an analogue of the fractional chain
rule of Christ and Weinstein, \cite{christw}.

\begin{proposition}[Fractional chain rule for the harmonic oscillator]\label{frac chain}
Let $H=-\tfrac{a}2\Delta + \tfrac{ib}2(x\cdot\nabla + \nabla\cdot x) + \tfrac{c}2 |x|^2$
with $a,c>0$ and normalized so that $ac-b^2=1$. Given $0\leq\gamma\leq\frac12$ and $1<p<\infty$,
then
$$
\bigl\| H^\gamma F(u) \bigr\|_p \lesssim \bigl\| H^\gamma u \bigr\|_{p_1} \, \bigl\| F'(u) \bigr\|_{p_2}
$$
provided $p_1,\ p_2\in(1,\infty)$ obey $\frac1{p_1}+\frac1{p_2}=\frac1p$.  The implicit constant depends continuously
on $a$, $b$, and $c$.  Here, $F(u)=|u|^{\frac4{n-2}}u$.
\end{proposition}

\begin{proof}
By Lemma~\ref{L:same norm}, the fractional chain rule from \cite{christw}, and H\"older's inequality, we have
\begin{align*}
\bigl\| H^\gamma F(u) \bigr\|_p
&\simeq \bigl\| (-\Delta)^\gamma F(u) \bigr\|_p + \bigl\| |x|^{2\gamma} F(u) \bigr\|_p \\
&\lesssim \bigl\| (-\Delta)^\gamma u \bigr\|_{p_1} \bigl\| F'(u) \bigr\|_{p_2}
    + \bigl\| |x|^{2\gamma} u \bigr\|_{p_1} \bigl\| F'(u) \bigr\|_{p_2} \\
&\lesssim \bigl\| H^\gamma u \bigr\|_{p_1} \bigl\| F'(u) \bigr\|_{p_2}.
\end{align*}
The last line is another application of Lemma~\ref{L:same norm}.
\end{proof}

\begin{corollary}\label{2.10}
Fix $n\geq 5$ and a time interval $I$. Then
$$
\Bigl\| \int_{t_0}^t U(t-s) F(u(s)) \,ds \Bigr\|_{W(I)} \lesssim \| u \|_{W(I)}^{\frac{n+2}{n-2}}.
$$
for each $t_0\in I$. In the case of the repulsive potential, the implicit constant depends continuously
on $T$ obeying $I\subseteq[-T,T]$.  In the case of a confining potential, it depends continuously on $|I|$.
In the free case, it is universal.
\end{corollary}

\begin{proof}
By \eqref{foschi est}, it suffices to show that
$$
\|F(u)\|_{N(I)}\lesssim \| u \|_{W(I)}^{\frac{n+2}{n-2}},
$$
which follows from Proposition~\ref{frac chain}, H\"older's inequality, and Lemma~\ref{sob em}.

The dependence of the constant in the statement of Corollary~\ref{2.10} is dictated by the Strichartz inequality in the confining case and by
\eqref{remark2.9} in the repulsive case.
\end{proof}

%%%%%%%%%%%%%%%%%%%%%%%%%%%%%%%%%%%%%%%%%%%%%%%%%%%%%%%%%%%%%%%%%%%%%%%%%%%%%%%%%%%%%%%%%%%
%
%
%                                   Section
%
%
%%%%%%%%%%%%%%%%%%%%%%%%%%%%%%%%%%%%%%%%%%%%%%%%%%%%%%%%%%%%%%%%%%%%%%%%%%%%%%%%%%%%%%%%%%%

\section{Monotonicity formulae}
In this section, we will prove two useful monotonicity estimates: a local mass conservation estimate and a Morawetz inequality.
The local mass conservation estimate will be used to control the flow of mass through a region of space.
As usual, the Morawetz inequality will be used to prevent concentration.

\subsection{Local mass conservation}
We recall a local mass conservation law that has appeared for instance in \cite{borg:scatter}, \cite{grillakis:scatter}, and
\cite{tao: gwp radial}. Let $\chi$ be a bump function supported on the ball $B(0,2)$ that equals $1$ on the ball $B(0,1)$.
We define
$$
\Mass(u(t), B(x_0,R))=\int \chi^2\Bigl(\frac{x-x_0}R\Bigr)|u(t,x)|^2 dx.
$$
Differentiating with respect to time we obtain
\begin{align*}
\partial_t \Mass(u(t), B(x_0,R))
&=2 \Re\int \chi^2\Bigl(\frac{x-x_0}R\Bigr) \bar u u_t(t,x) dx\\
&=-2 \Im \int\chi^2\Bigl(\frac{x-x_0}R\Bigr)\bar u\Delta u(t,x) dx.
\end{align*}
An integration by parts yields
$$
\partial_t \Mass(u(t), B(x_0,R))
=\frac 4R\int \chi\Bigl(\frac{x-x_0}R\Bigr)\nabla \chi \Bigl(\frac{x-x_0}R\Bigr)\Im (\bar u\nabla u)(x,t) dx,
$$
and so, by H\"older's inequality, we get
\begin{align}\label{rate mass change}
\partial_t \Mass(u(t), B(x_0,R))\lesssim \frac 1R \|\nabla u(t)\|_2 \Bigl(\Mass(u(t), B(x_0,R))\Bigr)^{\frac 12}.
\end{align}
On the other hand, by H\"older and Sobolev embedding we have
\begin{align}\label{mass in ball}
\Mass(u(t), B(x_0,R)) &\lesssim R^2\|\nabla u(t)\|_2^2,
\end{align}
which controls the concentration of mass in small volumes.

\subsection{A Morawetz inequality}\label{SS:morawetz}
Like other treatments of energy-critical problems, our analysis rests on an inequality that captures the defocusing sign of the nonlinearity.
The exact form this inequality takes varies from paper to paper, but in all cases it is based on the following formula for the derivative of the
momentum density:
\begin{equation}\label{momentum id}
\partial_t \Im(u_k \bar u)
=\tfrac14 \partial_k\Delta(|u|^2)-\Re\partial_j(\bar u_ku_j)-V_k|u|^2-\tfrac2n \partial_k|u|^{\dnnt}.
\end{equation}
Indeed, this will lead us to the following \emph{a priori} information:

\begin{proposition}\label{morawetz}
Let $u$ be a solution to \eqref{equation} on a spacetime slab $I\times\R^n$.  Then, for any $K\ge 1$, we have

\noindent 1. If $V\equiv0$, then
$$
\int_I\int_{|x|\le K|I|^{\frac12}}\frac{|u(t,x)|^{\dnnt}}{|x|}\,dx\, dt  \lesssim K|I|^{\frac 12} E(u).
$$

\noindent 2. If $V(x)=-\frac 12 |x|^2$, then
$$
\int_I\int_{|x|\le K|I|^{\frac12}}\frac{|u(t,x)|^{\dnnt}}{|x|}\,dx\, dt
\lesssim K|I|^{\frac 12} \Bigl(\|u\|_{L_t^{\infty}\dot H_x^1(I\times\R^n)}^2+ \|u\|_{L_t^{\infty}L_x^{\dnnt}(I\times\R^n)}^{\dnnt}\Bigr).
$$

\noindent 3. If $V(x)=\frac 12 |x|^2$, then
$$
\int_I\int_{|x|\le K|I|^{\frac12}}\frac{|u(t,x)|^{\dnnt}}{|x|}\,dx\,dt
\lesssim K\bigl( |I|^{\frac12} + |I| \bigr) E(u).
$$
\end{proposition}

\begin{proof}
For $V\equiv 0$, this result appears in \cite{borg:scatter,grillakis:scatter,tao:  gwp radial}; we will mimic the proof in \cite{tao: gwp radial}.

Let $\psi(x)$ be a radial nondecreasing (in radius) function obeying
\begin{equation*}
\psi(x)=
\begin{cases}
|x| &\text{if } |x|\leq 1\\
\tfrac32 &\text{if } |x|\geq 2,
\end{cases}
\end{equation*}
which is smooth except at the origin.  Given $R>0$, we define $a(x):=R\, \psi\bigl(\frac xR\bigr)$.

Integrating \eqref{momentum id} against the gradient of $a$ gives
\begin{align*}
\partial_t\int a_k(x)&\Im(u_k\bar u)(t,x)\,dx \\
&=\frac 14\int(-\Delta\Delta a)(x)|u(t,x)|^2\, dx + \int a_{jk}(x)\Re(\bar u_k u_j)(t,x)\,dx\\
&\quad -\int a_k(x)V_k(x)|u(t,x)|^2\,dx + \frac 2n\int \Delta a(x)|u(t,x)|^{\dnnt}\,dx.
\end{align*}
A few computations show that for $|x|\le R$ we have
\begin{align*}
a_k(x)&=\frac{x_k}{|x|}, \ \  \Delta a(x)\geq \frac{n-1}{|x|}> 0,\\
-\Delta\Delta a&>0 \quad \text{as a distribution,}
\end{align*}
and the matrix $a_{jk}(x)$ is positive definite.  Hence,
\begin{align*}
\partial_t\int&  a_k(x) \Im(u_k\bar u)  (t,x)\,dx\\
&\gtrsim  \int_{ R\leq |x|\leq 2R}(-\Delta\Delta a)(x)|u(t,x)|^2\, dx + \int_{ R\leq |x|\leq 2R} a_{jk}(x)\Re(\bar u_k u_j)(t,x)\,dx\\
&\quad -\int_{|x|\leq 2R} a_k(x)V_k(x)|u(t,x)|^2\,dx + \int_{ R\leq |x|\leq 2R} \Delta a(x)|u(t,x)|^{\dnnt}\,dx \\
&\quad + \int_{|x|\leq R} \frac{|u(t,x)|^{\dnnt}}{|x|}\,dx.
\end{align*}
In the region $R\le |x|\le 2R$, we have the rough estimates
\begin{align*}
|a_k(x)| \lesssim 1, \ |a_{jk}(x)|&\lesssim\frac 1R, \ \text{and}\ |-\Delta\Delta a(x)|\lesssim \frac 1{R^3}.
\end{align*}
Using these together with \eqref{mass in ball}, we obtain
\begin{gather*}
\biggl| \int_{ R\leq |x|\leq 2R}(-\Delta\Delta a)(x)|u(t,x)|^2\, dx \biggr|
    \lesssim R^{-3} \int_{|x|\leq 2R}|u(t,x)|^2 \lesssim R^{-1} \|\nabla u(t)\|_2^2, \\
\biggl| \int_{ R\leq |x|\leq 2R} a_{jk}(x)\Re(\bar u_k u_j)(t,x)\,dx \biggr|
    \lesssim R^{-1} \|\nabla u(t)\|_2^2, \\
\intertext{and}
\biggl| \int_{ R\leq |x|\leq 2R} \Delta a(x)|u(t,x)|^{\dnnt}\,dx \biggr|
    \lesssim R^{-1} \|u(t)\|_\dnnt^{\dnnt}.
\end{gather*}
Also, by Cauchy--Schwarz and \eqref{mass in ball},
\begin{align*}
\biggl| \int  a_k(x) \Im(u_k\bar u)  (t,x)\,dx\biggr|  \lesssim \biggl| \int_{|x|\leq 2R} |\nabla u(t,x)|\,|u(t,x)| \,dx\biggr|
    \lesssim R \|\nabla u(t)\|^2_2.
\end{align*}
This leaves us to estimate the term involving the potential $V$.  As $a(x)$ is a non-decreasing function of radius,
\begin{align*}
-\int_{|x|\leq 2R} a_k(x)V_k(x)|u(t,x)|^2\,dx \geq 0,
\end{align*}
for any repulsive potential, that is, obeying $x\cdot\nabla V(x)\leq 0$.  In particular, this covers $V(x)=-\frac12 |x|^2$ and
$V\equiv 0$.  For $V(x)=\frac12 |x|^2$, using \eqref{mass in ball} we obtain
\begin{align*}
\biggl| \int_{|x|\leq 2R} a_k(x)V_k(x)|u(t,x)|^2\,dx \biggr|
&\lesssim \int_{|x|\le 1}|u(t,x)|^2\,dx\,dt + \int_{|x|\ge 1}|x|^2|u(t,x)|^2\,dx\,dt\\
&\lesssim \|\nabla u(t)\|_2^2 + \|xu(t)\|_2^2.
\end{align*}
The proposition follows from these estimates and the Fundamental Theorem of Calculus, by choosing $R=K|I|^{1/2}$.
\end{proof}

%%%%%%%%%%%%%%%%%%%%%%%%%%%%%%%%%%%%%%%%%%%%%%%%%%%%%%%%%%%%%%%%%%%%%%%%%%%%%%%%%%%%%%%%%%%
%
%
%                                   Section
%
%
%%%%%%%%%%%%%%%%%%%%%%%%%%%%%%%%%%%%%%%%%%%%%%%%%%%%%%%%%%%%%%%%%%%%%%%%%%%%%%%%%%%%%%%%%%%

\section{Local theory}

The local theory for the energy-critical NLS ($V\equiv 0$) was mostly worked out by Cazenave and Weissler, \cite{cwI}.
They constructed local-in-time solutions for arbitrary initial data in $\dot H^1(\R^n)$ and global solutions for small energy data.
However, as for any critical equation, the time of existence of the local solution depends on the profile of the initial data and not
simply on its $\dot H^1_x$-norm; this is the reason why global existence does not follow immediately from the conservation of energy
and the usual iterative argument.  They also proved uniqueness of these solutions in certain Strichartz spaces; moreover, the solution
was shown to depend continuously (in these spaces) on the initial data in the energy space $\dot H^1(\R^n)$.
A later argument of Cazenave, \cite{cazenave:book}, also demonstrates that the uniqueness is in fact unconditional in the category of
strong solutions, that is, solutions belonging to $C_t^0\dot H^1_x$. Finally, Tao and Visan, \cite{TV}, showed that the map from initial
data to strong solutions is uniformly continuous in energy-critical spaces; see also \cite{ckstt:gwp, RV} for a proof in dimensions $n=3,4$.

In this section we develop a local theory for the energy-critical NLS with harmonic potential. The first step is to construct local-in-time solutions.

\begin{proposition}[Local well-posedness]\label{lwp vn0}
Let $V(x)=\pm \frac 12 |x|^2$, $u_0\in\Sigma$, and $I$ be a compact time interval that contains $0$ such that
\begin{align}\label{lwp ass}
\|U(t)\nabla u_0\|_{L_t^{\frac{2(n+2)}{n-2}}L_x^{\frac{2n(n+2)}{n^2+4}}(\ir)}\le \eta,
\end{align}
for a sufficiently small absolute constant $\eta>0$.  In the confining case, we also make an assumption on the size of the interval, say $|I|<1$.
Then, there exists a unique strong solution to \eqref{equation} on $\ir$ such that
\begin{align}\label{stric<energy}
\|H(-t)^{\frac 12} u(t)\|_{ S^0(I)}\le C(\|u_0\|_\Sigma).
\end{align}
(Compare \eqref{stric<energy} with \eqref{H(-t) control free}, \eqref{H(-t) control conf}, and \eqref{H(-t) control repuls}.)
\end{proposition}

\begin{proof} The proof of Proposition~\ref{lwp vn0} is standard and based on contraction mapping arguments.  We define the solution map to be
$$
[\Phi(u)](t):=U(t)u_0 - i\int_0^t U(t-s)F(u(s))\, ds.
$$
Thus, applying the momentum and position operators we find
$$
P(-t)[\Phi(u)](t)=U(t)i\nabla u_0 - i\int_0^t U(t-s)P(-s)F(u(s))\, ds
$$
and
$$
X(-t)[\Phi(u)](t)=U(t)xu_0 - i\int_0^t U(t-s)X(-s)F(u(s))\, ds.
$$
By Strichartz, Lemma~\ref{lemma jt ht}, Lemma~\ref{lemma pt sobolev}, and \eqref{lwp ass}, we estimate
\begin{align*}
\|P(-t)\Phi(u)&\|_{L_t^{\frac{2(n+2)}{n-2}}L_x^{\frac{2n(n+2)}{n^2+4}}}\\
&\lesssim  \eta + \|P(-t)u\|_{L_t^{\frac{2(n+2)}{n-2}}L_x^{\frac{2n(n+2)}{n^2+4}}}^{\frac{n+2}{n-2}} \\
\|X(-t)\Phi(u)&\|_{L_t^{\frac{2(n+2)}{n-2}}L_x^{\frac{2n(n+2)}{n^2+4}}}\\
&\lesssim \|xu_0\|_2 + \|X(-t) u\|_{L_t^{\frac{2(n+2)}{n-2}}L_x^{\frac{2n(n+2)}{n^2+4}}}
        \|P(-t)u\|_{L_t^{\frac{2(n+2)}{n-2}}L_x^{\frac{2n(n+2)}{n^2+4}}}^{{\frac4{n-2}}}\\
        \|\Phi(u)&\|_{L_t^{\frac{2(n+2)}{n-2}}L_x^{\frac{2n(n+2)}{n^2+4}}}\\
&\lesssim \|u_0\|_2 + \|u\|_{L_t^{\frac{2(n+2)}{n-2}}L_x^{\frac{2n(n+2)}{n^2+4}}}
        \|P(-t)u\|_{L_t^{\frac{2(n+2)}{n-2}}L_x^{\frac{2n(n+2)}{n^2+4}}}^{{\frac4{n-2}}}
\end{align*}
where all spacetime norms are on $\ir$.  Similarly,
\begin{align*}
\|&\Phi( u)  -\Phi(v)\|_{L_t^{\frac{2(n+2)}{n-2}}L_x^{\frac{2n(n+2)}{n^2+4}}}\\
&\lesssim \Bigl( \|P(-t)u\|_{L_t^{\frac{2(n+2)}{n-2}}L_x^{\frac{2n(n+2)}{n^2+4}}}^{{\frac4{n-2}}}
        + \|P(-t)v\|_{L_t^{\frac{2(n+2)}{n-2}}L_x^{\frac{2n(n+2)}{n^2+4}}}^{{\frac4{n-2}}}\Bigr)
    \|u-v\|_{L_t^{\frac{2(n+2)}{n-2}}L_x^{\frac{2n(n+2)}{n^2+4}}}.
\end{align*}
It is thus easy to see that $\Phi$ maps
\begin{align*}
\mathcal B=\{u;\ & \|P(-t) u\|_{L_t^{\frac{2(n+2)}{n-2}}L_x^{\frac{2n(n+2)}{n^2+4}}(\ir)}\le 2\eta, \\
                 & \|X(-t) u\|_{L_t^{\frac{2(n+2)}{n-2}}L_x^{\frac{2n(n+2)}{n^2+4}}(\ir)}\le 2C \|xu_0\|_2,\\
                 & \|u\|_{L_t^{\frac{2(n+2)}{n-2}}L_x^{\frac{2n(n+2)}{n^2+4}}(\ir)}\le 2C\|u_0\|_2\}
\end{align*}
to itself and is a contraction in the $L_t^{\frac{2(n+2)}{n-2}}L_x^{\frac{2n(n+2)}{n^2+4}}$ norm,
provided $\eta$ is chosen sufficiently small depending only on the Strichartz constant $C$ and on the Sobolev embedding constant
in Lemma~\ref{lemma pt sobolev}.  The contraction mapping theorem then implies the existence of a unique solution to \eqref{equation} on $I$.
As $H(t)=\tfrac 12 P(t)^2 + \tfrac12 X(t)^2$, the claim \eqref{stric<energy} also follows.
\end{proof}
\begin{remark}
An immediate consequence of Proposition~\ref{lwp vn0} and Strichartz inequality is global well-posedness for \eqref{equation}
in the repulsive case for initial data small in $\dot H^1(\R^n)$.  Scattering for small $\dot H^1(\R^n)$ initial data in this case
is a consequence of \eqref{stric<energy}, see Section~\ref{scattering section}.
\end{remark}
\begin{remark}
It is not hard to see that in the confining case global well-posedness also holds under a small energy assumption.
However, as the global solution is constructed by stacking up intervals,  this does not imply scattering.  Indeed, we believe
that scattering (to the harmonic oscillator) fails.  Moreover, it is unlikely that a smallness assumption on the $\dot H^1_x$-norm of
the initial data (combined with finite $\Sigma$-norm) is sufficient to construct a global solution by such a simple argument.
\end{remark}
\begin{remark}
For large energy initial data, \eqref{lwp ass} is satisfied for small time intervals $I$ as can be seen from the Strichartz inequality and
the Dominated Convergence Theorem.  However, this implies that the time of existence depends on the profile of the initial data
rather than its energy alone.
\end{remark}

Next we show that in the presence of a quadratic potential, the Strichartz norms of a solution can be bounded in terms of its $Z$-norm.

\begin{lemma}\label{stric<Z}
Let $I$ be a compact time interval containing $0$.  Suppose $V=\pm\frac12|x|^2$ and $u$ is
a strong solution to \eqref{equation} on $\ir$ with $ \|u\|_{Z(I)}<L$.  Then
$$
\| H(-t)^{\frac 12} u\|_{ S^0(I)}\leq C(L, |I|, \|u_0\|_\Sigma).
$$
\end{lemma}
\begin{proof} Let $\eta$ be a small constant to be specified later.  Subdivide $I$ into $N=N(L,|I|)$ smaller intervals $I_j=[t_{j}, t_{j+1}]$
such that on each such interval,
$$
\|u\|_{Z(I_j)}\sim \eta.
$$
In the confining case, we assume in addition that $|I_j|\le 1$.

Applying the Strichartz estimate in the Duhamel formula and using Proposition~\ref{frac chain} gives
\begin{align*}
\|H(-t)^{\frac 12}u\|_{S^0(I_j)}
&\lesssim\|H(-t_j)^{\frac 12}u(t_j)\|_2 + \|u\|_{Z(I_j)}^{{\frac4{n-2}}}\|H(-t)^{\frac 12}u\|_{S^0(I_j)}\\
&\lesssim \|H(-t_j)^{\frac 12}u(t_j)\|_2 + \eta^{{\frac4{n-2}}}\|H(-t)^{\frac 12}u\|_{S^0(I_j)},
\end{align*}
which by choosing $\eta$ sufficiently small depending only on the Strichartz constant, implies
\begin{align*}
\|H(-t)^{\frac 12}u\|_{S^0(I_j)} \lesssim \|H(-t_j)^{\frac 12}u(t_j)\|_2 \leq C( \|u_0\|_\Sigma).
\end{align*}
In the confining case, the last inequality follows from the conservation of energy, see \eqref{H(-t) control conf}; in the repulsive case,
see \eqref{H(-t) control repuls}.  The claim follows by adding these estimates over all subintervals $I_j$.
\end{proof}

Next, we establish a standard blowup criterion for solutions to \eqref{equation} in the presence of a potential.
In the free case, this is due to Cazenave and Weissler, \cite{cwI}.

\begin{lemma}[Blowup criterion]\label{blow}
Let $V=\pm\frac12|x|^2$, $u_0\in \Sigma$, and let $u$ be a strong solution to \eqref{equation} on the slab $[0, T_0)\times\R^n$ such that
\begin{align}\label{norm finite}
\|u\|_{Z([0, T_0))} < \infty.
\end{align}
Then there exists $\delta>0$ such that the solution $u$ extends to a strong solution to \eqref{equation} on the slab
$[0, T_0+\delta]\times\R^n$.
\end{lemma}

In the contrapositive, this lemma asserts that if a solution cannot be continued strongly beyond a time $T_*$, then the
$Z$-norm must blow up at that time.  One can also establish that other scale-invariant norms
(except for those norms involving $L^\infty_t$) also blow up at this time, but we will not do so here.

\begin{proof}
For $t_0\in [0,T_0)$ and $t_0\leq t< T_0$ we have the Duhamel formula
\begin{align*}
U(t-t_0)u(t_0)=u(t)- i \int_{t_0}^t U(t-s) F(u(s))\, ds.
\end{align*}
Thus,
\begin{align*}
\|U(t-t_0)\nabla u(t_0)\|_{ S^0([t_0, T_0))}
&=\|P(-t) U(t-t_0) u(t_0)\|_{ S^0([t_0, T_0))} \\
&\lesssim_{T_0}\|H(-t)^{\frac 12} U(t-t_0) u(t_0)\|_{ S^0([t_0, T_0))} \\
&\lesssim_{T_0} \|H(-t)^{\frac 12} u\|_{ S^0([t_0, T_0))} + \|H(-t)^{\frac 12} u\|_{ S^0([t_0, T_0))}^{\frac{n+2}{n-2}},
\end{align*}
where the first inequality follows from Lemma~\ref{L:same norm} and the second from Strichartz, Proposition~\ref{frac chain}, and \eqref{emb}.
As by \eqref{norm finite} and Lemma~\ref{stric<Z},
$$
\|H(-t)^{\frac 12} u\|_{ S^0([0, T_0))}<\infty,
$$
the Dominated Convergence Theorem implies that we may choose $t_0$ close enough to $T_0$ so that
\begin{align*}
\|U(t-t_0)\nabla u(t_0)\|_{L_t^{\frac{2(n+2)}{n-2}}L_x^{\frac{2n(n+2)}{n^2+4}}([t_0, T_0)\times\R^n)}\leq\tfrac 12 \eta,
\end{align*}
where $\eta$ is as in Proposition~\ref{lwp vn0}.  From the control we have on $\|\nabla u(t)\|_2$ (see subsection~\ref{energy control ss}), Strichartz
inequality, and the Dominated Convergence Theorem, one may choose $\delta>0$ such that
\begin{align*}
\|U(t-t_0)\nabla u(t_0)\|_{L_t^{\frac{2(n+2)}{n-2}}L_x^{\frac{2n(n+2)}{n^2+4}}([T_0, T_0+\delta]\times\R^n)}\leq \tfrac 12 \eta.
\end{align*}
Thus, we may apply Proposition~\ref{lwp vn0} on the interval $[t_0, T_0+\delta]$ to complete the proof.
\end{proof}

%%%%%%%%%%%%%%%%%%%%%%%%%%%%%%%%%%%%%%%%%%%%%%%%%%%%%%%%%%%%%%%%%%%%%%%%%%%%%%%%%%%%%%%%%%%
%
%
%                                   Section
%
%
%%%%%%%%%%%%%%%%%%%%%%%%%%%%%%%%%%%%%%%%%%%%%%%%%%%%%%%%%%%%%%%%%%%%%%%%%%%%%%%%%%%%%%%%%%%

\section{A perturbative result}\label{S:pert}

In this section we show that the nonlinear solution cannot be large without the linear solution being large.  This is an analogue
of Lemma~3.2 in \cite{tao: gwp radial}, which deals with $V\equiv 0$.  In low dimensions, $n=3,4,5$, the proof is essentially identical.
However, in higher dimensions, the proof in \cite{tao: gwp radial} is rather complicated; the treatment here is significantly simpler.
The main simplification arises from working in exotic Strichartz spaces with a fractional number of derivatives (but critical scaling);
this is reminiscent of the approaches used in \cite{nak, TV}.  Pushing this through requires the exotic Strichartz estimates of Foschi,
\eqref{foschi est}.

\begin{lemma}[Perturbation lemma]\label{pert lemma}
Let $u$ be a solution to \eqref{equation} on $I=[t_1,t_2]$ such that
\begin{align}\label{pert ass}
\tfrac 12 \eta \leq \|u\|_{W(I)} \leq \eta,
\end{align}
where $\eta$ is a sufficiently small constant depending on the norm of the initial data and, in the presence of a potential, on $T$ defined by
$I\subseteq [-T,T]$.  Then
\begin{align*}
\|u_k\|_{W(I)}\geq \tfrac14 \eta,
\end{align*}
where $u_k(t)=U(t-t_k)u(t_k)$ and $k=1,2$.
\end{lemma}

\begin{proof}
We will present the proof for $k=1$; the other case is basically identical.  We begin with dimensions $n=3,4,5$.

Let $V(x)=\pm\tfrac 12 |x|^2$.  By \eqref{emb},
$$
\|u\|_{Z(I)}\lesssim_T \|u\|_{W(I)}\lesssim_T \eta,
$$
and hence, by Lemma~\ref{stric<Z},
$$
\| H(-t)^{\frac12} u \|_{S^0(I)}\leq C(T,\|u_0\|_\Sigma).
$$
Applying \eqref{emb} and Strichartz, we obtain
\begin{align*}
\|u-u_1\|_{W(I)}
&\lesssim_T \| H(-t)^{\frac12} (u-u_1) \|_{ S^0(I)} \\
&\lesssim_T \| H(-t)^{\frac12} F(u) \|_{L_t^2 L_x^\frac{2n}{n+2}(I\times\R^n)} \\
&\lesssim_T \| H(-t)^{\frac12} u \|_{ S^0(I)}  \| u \|_{Z(I)}^{\frac{4}{n-2}} \\
&\leq C(T,\|u_0\|_\Sigma) \eta^{\frac{4}{n-2}}.
\end{align*}
As $\tfrac{4}{n-2}>1$ in dimensions $n=3,4,5$, the result follows from \eqref{pert ass} and the triangle inequality,
provided $\eta$ is chosen sufficiently small depending on $T$ and $\| u_0\|_\Sigma$.

Of course, in the free case the same argument yields
$$
\|u-u_1\|_{W(I)} \leq C(\|\nabla u_0\|_2) \eta^{\frac{4}{n-2}},
$$
so the claim follows as before choosing $\eta$ sufficiently small depending only on $\| \nabla u_0\|_2$.

We now turn to dimensions $n\geq 5$.  By Duhamel's formula, Corollary~\ref{2.10}, and \eqref{pert ass},
\begin{align*}
\|u-u_1\|_{W(I)}
&\lesssim \| u \|_{W(I)}^{\frac{n+2}{n-2}}
\lesssim \eta^{\frac{n+2}{n-2}},
\end{align*}
where the implicit constants depend on $T$ in the presence of a potential.  Once again, the result follows from \eqref{pert ass}
and the triangle inequality, provided $\eta$ is chosen small enough.
\end{proof}

%%%%%%%%%%%%%%%%%%%%%%%%%%%%%%%%%%%%%%%%%%%%%%%%%%%%%%%%%%%%%%%%%%%%%%%%%%%%%%%%%%%%%%%%%%%
%
%
%                                   Section
%
%
%%%%%%%%%%%%%%%%%%%%%%%%%%%%%%%%%%%%%%%%%%%%%%%%%%%%%%%%%%%%%%%%%%%%%%%%%%%%%%%%%%%%%%%%%%%

\section{Global well-posedness}\label{gwp section}

In this section, we prove global well-posedness for large data.  This settles the theorems stated in the introduction,
with the exception of scattering, which is treated in Section~\ref{scattering section}.
We follow \cite{tao: gwp radial} almost verbatim, incorporating a few simplifications that result from the improved perturbation theory
described in Section~\ref{S:pert}.

\subsection{The setup}\label{SS:6.1}
In the confining case, it suffices to show that the Cauchy problem \eqref{equation} with initial data in $\Sigma$ admits a unique
solution on a small time interval containing zero, say $I=[-1,1]$.  Global well-posedness follows from the conservation of energy and
the usual iterative argument.

In the free case, we show that given initial data $u_0\in \dot H^1(\R^n)$ and a compact time interval $I=[-T, T]$,
the Cauchy problem \eqref{equation} admits a unique solution $u$ on $I$ and moreover,
\begin{align}\label{goal free}
\|u\|_{Z(I)}\leq C(E(u_0)).
\end{align}
Global well-posedness is then a consequence of this fact and Lemma~\ref{blow}.  Scattering is an easy consequence of \eqref{goal free};
see, for example, \cite{cazenave:book, TV}.

In the repulsive case, given initial data $u_0\in \Sigma$ and a compact time interval $I=[-T,T]$, the Cauchy problem
\eqref{equation} will also be shown to admit a unique solution $u$ on $I$.  In this case, we will establish
\begin{align}\label{goal repuls}
\|u\|_{Z(I)}\leq C(T,\|u_0\|_{\Sigma}).
\end{align}
Global well-posedness follows again from this and Lemma~\ref{blow}; however, \eqref{goal repuls} alone is not strong enough
to imply scattering.  See Section~\ref{scattering section} for a proof of scattering in the repulsive case.

Thus, let
\begin{equation*}
[t_-,t_+]=
\begin{cases}
[-1,1] & \quad \text{if} \quad  V(x)=\tfrac 12 |x|^2\\
[-T,T] & \quad \text{if} \quad  V\equiv 0 \text{ or } V(x)=-\tfrac 12 |x|^2.
\end{cases}
\end{equation*}
We will prove that the Cauchy problem \eqref{equation} admits a unique solution $u$ on $[t_-,t_+]$.  By Lemma~\ref{blow}, it
suffices to \emph{a priori} assume that the solution already exists on $[t_-,t_+]$ and prove that it obeys
\begin{align}\label{goal}
\|u\|_{Z([t_-,t_+])}\le
\begin{cases}
C(E(u_0)) & \quad \text{if} \quad  V\equiv 0 \text{ or } V(x)=\tfrac 12 |x|^2 \\
C(T, \|u_0\|_\Sigma) & \quad \text{if} \quad  V(x)=-\tfrac 12 |x|^2.
\end{cases}
\end{align}
In view of this goal, we adopt the following

\begin{convent}
In this section, all implicit constants are permitted to depend on the dimension $n$, the energy in the free and confining cases,
and on $T$ and $\|u_0\|_\Sigma$ in the repulsive case.
\end{convent}

To this end, we divide $[t_-,t_+]$ into $J$ subintervals $I_j=[t_j,t_{j+1}]$ such that
\begin{align}\label{small W}
\frac 12 \eta \leq \|u\|_{W(I_j)}\leq \eta,
\end{align}
where $\eta$ is a small constant depending on the dimension $n$, the energy in the free and confining cases, and on $\|u_0\|_{\Sigma}$ and $T$
in the repulsive case.  Indeed, in this section this is what we will always mean by the phrase `sufficiently small'.
Note that \eqref{goal} is equivalent to estimating $J$, and this is what we will do.

\begin{remark}
By Sobolev embedding (cf. \eqref{emb}) and \eqref{small W},
\begin{align}\label{small Z}
\|u\|_{Z(I_j)}\lesssim \eta.
\end{align}
Combining this with Proposition~\ref{stric<Z} gives
\begin{align}\label{H OK}
\|H(-t)^{\frac 12}u\|_{S^0(I_j)}\lesssim 1.
\end{align}
\end{remark}

Now let $u_\pm:= U(t-t_\pm)u(t_\pm)$.  Note that by Sobolev embedding, Strichartz, and \emph{a priori} energy
control (see \eqref{H(-t) control free}, \eqref{H(-t) control conf}, and \eqref{H(-t) control repuls}),
\begin{equation}\label{bdd}
\begin{aligned}
\|u_\pm\|_{Z([t_-,t_+])}
&\lesssim \|u_\pm\|_{W([t_-,t_+])}
\lesssim \|H(-t)^{\frac 12} u_\pm\|_{ S^0([t_-,t_+])} \\
&\lesssim \|H(-t_\pm)^{\frac 12} u(t_\pm)\|_2
\lesssim 1.
\end{aligned}
\end{equation}

As in \cite{tao: gwp radial} we make the following

\begin{definition}
We call $I_j$ exceptional if
$$
\|u_\pm\|_{W(I_j)} >\eta^{C_1}
$$
for at least one sign $\pm$.  Otherwise, we call $I_j$ unexceptional.
\end{definition}

The choice of $C_1$ in this definition will need to satisfy several constraints as we proceed.  Eventually we will find that we may choose
$C_1= 15 n^2$, for example.

\begin{remark}\label{count ex}
The estimate \eqref{bdd} gives an upper bound on the number of exceptional intervals, $O(\eta^{-\frac{2(n+2)}{n-2}C_1})$.  Moreover, if there were no
unexceptional intervals, the claim \eqref{goal} would follow from this bound and \eqref{small Z}.  As a consequence,
we need only count the number of unexceptional intervals.
\end{remark}

By this remark, we may assume that there exist unexceptional intervals $I_j$.  As in \cite{tao: gwp radial} and in its predecessor, \cite{borg:scatter},
the first step in counting such intervals is to prove that each contains a bubble of mass concentration.

\subsection{A concentration result}

\begin{proposition}[Existence of a bubble]\label{bubble}
Let $I_j$ be an unexceptional interval.  Then there exists $x_j\in\R^n$ such that
$$
\Mass\bigl( u(t), B(x_j, \eta^{- C} |I_j|^{\frac12} ) \bigr) \gtrsim \eta^{C} |I_j|
$$
for all $t\in I_j$.  Here we may take $C=4 n^2$.
\end{proposition}

\begin{proof}
Fix $I_j=[t_j,t_{j+1}]$ unexceptional and let $t_*=\tfrac12(t_j+t_{j+1})$ be its mid-point.  We
will use the abbreviations $I_j^1=[t_j,t_*]$ and $I_j^2=[t_*,t_{j+1}]$ for the two halves of $I_j$.

As $\|u\|_{W(I_j)} \geq \tfrac 12 \eta$, we may invoke time-reversal symmetry to assume that $\|u\|_{W(I_j^{2})} \geq \frac14 \eta$.
Thus by Lemma~\ref{pert lemma},
\begin{equation}\label{big free}
\bigl\| U(t-t_*) u(t_*) \bigr\|_{W(I_j^2)} \geq \tfrac18 \eta.
\end{equation}

By Duhamel's formula,
\begin{equation}\label{**}
\begin{aligned}
U(t-t_*) u(t_*)
&=  U(t-t_-) u(t_-) - i \int_{t_-}^{t_j}   U(t-s) F(u(s))\,ds \\
&\quad        - i \int_{t_j}^{t_*}   U(t-s) F(u(s))\,ds.
\end{aligned}
\end{equation}
As $u_-:=U(t-t_-)u(t_-)$ and $I_j$ is unexceptional,
$$
\|U(t-t_-) u(t_-)\|_{W(I_j^2)} = \|u_-\|_{W(I_j^2)} \leq \eta^{C_1}.
$$
Moreover, by Corollary~\ref{2.10},
$$
\biggl\| \int_{t_j}^{t_*}   U(t-s) F(u(s))\,ds \biggr\|_{W(I_j^2)} \lesssim \|u\|_{W(I_j^1)}^{(n+2)/(n-2)}
    \lesssim \eta^{(n+2)/(n-2)}.
$$
Thus by the triangle inequality,
$$
\biggl\| \int_{t_-}^{t_j} U(t-s) F(u(s))\,ds \biggr\|_{W(I_j^2)}
    \geq \tfrac1{100}\eta,
$$
provided $\eta$ is chosen sufficiently small.  Here, we are also using the fact that $C_1$ and $\frac{n+2}{n-2}$ are both larger than $1$.
Henceforth, we will write
$$
v(t) := \int_{t_-}^{t_j} U(t-s) F(u(s))\,ds.
$$
In particular, we may rewrite the previous inequality as
\begin{equation}\label{v lower bound}
\| v \|_{W(I_j^2)} \geq \tfrac1{100}\eta.
\end{equation}

Next, we obtain an upper bound on $v$.  By applying the triangle inequality in \eqref{**} and using Strichartz, \eqref{small Z}, and \eqref{H OK},
we estimate
\begin{align}
\bigl\| H(-t)^\frac12 v(t) \bigr\|_{ S^0(I_j^2)}
&\lesssim \bigl\| U(t-t_*) H(-t_*)^\frac12 u(t_*) \bigr\|_{ S^0(I_j^2)}
+ \bigl\| H(-t)^\frac12 u_-(t) \bigr\|_{ S^0(I_j^2)}  \notag \\
&\quad + \biggl\| \int_{t_j}^{t_*}  U(t-s) H(-s)^\frac12 F(u(s))\,ds \biggr\|_{ S^0(I_j^2)} \notag \\
&\lesssim \bigl\| H(-t_*)^\frac12 u(t_*) \bigr\|_{2} + \bigl\| H(-t_-)^\frac12 u(t_-) \bigr\|_{2} \notag \\
&\quad + \| H(-t)^\frac12 u(t) \|_{S^0(I_j^1)} \| u \|_{Z(I_j^1)}^{\frac4{n-2}} \notag \\
&\lesssim 1. \label{v in S1}
\end{align}
Using this and \eqref{v lower bound}, interpolation gives
\begin{equation}\label{v lower bound 2}
\| v \|_{Z(I_j^2)} \gtrsim \eta^{n/2}.
\end{equation}

The next ingredient is a regularity result for $v$.  This mirrors Lemma~3.4 in \cite{tao: gwp radial}, but we streamline the proof.

\begin{lemma}\label{L:boho}
Let us write $v_{av}$ for the convolution of $v$ with $r^{-n}\varphi\bigl(\frac{\cdot}{r}\bigr)$ where $r= \eta^{C_0} |I_j|^{\frac 12}$ and $\varphi$ is such that its Fourier transform is smooth, compactly supported,
and equal to one in a neighborhood of the origin.  Then
\begin{align*}
\|v-v_{av}\|_{Z(I_j^2)}\lesssim \eta^{\frac{n-2}{n+2}C_0}.
\end{align*}
Here we may choose $C_0=3n$.
\end{lemma}

\begin{proof}
The claim will follow from H\"older's inequality once we establish
\begin{align}\label{boho}
\|v-v_{av}\|_{L_t^\infty L_x^{\frac{2(n+2)}{n-2}}(I_j^2\times\R^n)}\lesssim \eta^{\frac{n-2}{n+2}C_0} |I_j|^{-\frac{n-2}{2(n+2)}}.
\end{align}

By the Mikhlin multiplier theorem and Lemma~\ref{L:same norm}, we estimate
\begin{align}\label{bohoho}
\|v-v_{av}\|_{L_t^\infty L_x^{\frac{2(n+2)}{n-2}}(I_j^2\times\R^n)}
&\lesssim r^{\frac{n-2}{n+2}} \bigl\| |\nabla|^{\frac{n-2}{n+2}} v\bigr\|_{L_t^\infty L_x^{\frac{2(n+2)}{n-2}}(I_j^2\times\R^n)} \notag\\
&\lesssim \eta^{\frac{n-2}{n+2}C_0} |I_j|^{\frac{n-2}{2(n+2)}}\bigl\| H(-t)^{\frac{n-2}{2(n+2)}} v(t)\bigr\|_{L_t^\infty L_x^{\frac{2(n+2)}{n-2}}(I_j^2\times\R^n)}.
\end{align}
To continue, we use the dispersive estimate \eqref{disp est} followed by Proposition~\ref{frac chain} and Lemmas~\ref{sob em} and \ref{L:same norm}:
\begin{align*}
\bigl\| H(-t)^{\frac{n-2}{2(n+2)}} &v\bigr\|_{L_t^\infty L_x^{\frac{2(n+2)}{n-2}}(I_j^2\times\R^n)}\\
&\lesssim \biggl\|\int_{t_-}^{t_j} |t-s|^{-\frac{2n}{n+2}}\bigl\| H(-s)^{\frac{n-2}{2(n+2)}} F(u(s))\bigr\|_{\frac{2(n+2)}{n+6}}\, ds\biggr\|_{L_t^\infty(I_j^2)}\\
&\lesssim |I_j|^{-\frac{n-2}{n+2}} \|u\|_{L_t^\infty L_x^{\frac{2n}{n-2}}}^{\frac4{n-2}} \bigl\| H(-t)^{\frac{n-2}{2(n+2)}} u(t)\bigr\|_{L_t^\infty L_x^{\frac{2n(n+2)}{(n-2)(n+4)}}}\\
&\lesssim |I_j|^{-\frac{n-2}{n+2}} \|\nabla u\|_{L_t^\infty L_x^2}^{\frac{4}{n-2}} \bigl\| H(-t)^{\frac12} u(t)\bigr\|_{L_t^\infty L_x^2}\\
&\lesssim |I_j|^{-\frac{n-2}{n+2}}\bigl\| H(-t)^{\frac12} u(t)\bigr\|_{L_t^\infty L_x^2}^{\frac{n+2}{n-2}},
\end{align*}
where all spacetime norms are on $I_j^2\times\R^n$.  Invoking \eqref{H(-t) control free}, \eqref{H(-t) control conf}, and \eqref{H(-t) control repuls}  we thus have
$$
\bigl\| H(-t)^{\frac{n-2}{2(n+2)}} v(t)\bigr\|_{L_t^\infty L_x^{\frac{2(n+2)}{n-2}}(I_j^2\times\R^n)} \lesssim |I_j|^{-\frac{n-2}{n+2}}.
$$
Combining this with \eqref{bohoho}, we derive \eqref{boho}.
\end{proof}

We now return to the proof of Proposition~\ref{bubble}.  By Lemma~\ref{L:boho} and \eqref{v lower bound 2},
\begin{align}\label{v av in Z}
\|v_{av}\|_{Z(I_j^2)}\gtrsim \eta^{\frac n2}-\eta^{\frac{n-2}{n+2}C_0}\gtrsim \eta^{\frac n2},
\end{align}
provided $C_0$ is chosen sufficiently large, for example, $C_0=3n$.
On the other hand, by H\"older's and Young's inequalities and then \eqref{v in S1},
\begin{align*}
\| v_{av} \|_{L_{t,x}^{\frac{2n}{n-2}}(I_j^2\times\R^n)}
&\lesssim |I_j|^{\frac{n-2}{2n}} \| v_{av} \|_{L^\infty_t L_{x}^{\frac{2n}{n-2}}(I_j^2\times\R^n)} \\
&\lesssim |I_j|^{\frac{n-2}{2n}} \| v \|_{L^\infty_t L_{x}^{\frac{2n}{n-2}}(I_j^2\times\R^n)} \\
&\lesssim |I_j|^{\frac{n-2}{2n}}.
\end{align*}
Noting that $Z = L_{t,x}^{\frac{2(n+2)}{n-2}}$ interpolates between $L_{t,x}^{\frac{2n}{n-2}}$ and
$L_{t,x}^\infty$ and then using the bound above together with \eqref{v av in Z}, we obtain
$$
\| v_{av} \|_{L_{t,x}^{\infty}(I_j^2\times\R^n)} \gtrsim \eta^{\frac{n(n+2)}4} |I_j|^{-\frac{n-2}4}.
$$
Thus there exists $(s_j,x_j)\in I_j^2\times \R^n$ such that
$$
v_{av}(s_j,x_j) \gtrsim \eta^{\frac{n(n+2)}4} |I_j|^{-\frac{n-2}4},
$$
and so, by Cauchy--Schwarz and the definition of $v_{av}$,
$$
\Mass\bigl( v(s_j), B(x_j, \eta^{C_0} |I_j|^{\frac12})\bigr)  \gtrsim \eta^{C} |I_j|
$$
where $C\geq \frac{n(n+2)}{2} + n C_0$, for example, we may take $C_0 = 3 n$ and $C = 4 n^2$.
By \eqref{rate mass change} this mass remains on a somewhat larger ball throughout $I_j$; more
precisely, taking $R= \eta^{-C} |I_j|^{\frac12} $ in \eqref{rate mass change} and choosing $\eta$ sufficiently small gives
\begin{equation}\label{v bubble}
\Mass\bigl( v(t), B(x_j, \eta^{- C} |I_j|^{\frac12})\bigr)  \gtrsim \eta^{C} |I_j|,
\end{equation}
for all $t\in I_j$.

The last step is to show that this mass concentration holds for $u$, not merely $v$.  Recall that
$u(t_j)=u_-(t_j) - i v(t_j)$.  We first show mass concentration for $u$ at time $t_j$.

As $I_j$ is unexceptional, there is a $\tau_j\in I_j$ so that
$$
\| u_{-}(\tau_j) \|_{\frac{2(n+2)}{n-2}} \lesssim \eta^{C_1} |I_j|^{-\frac{n-2}{2(n+2)}}
$$
and so by H\"older's inequality,
\begin{align*}
\Mass\bigl( u_{-}(\tau_j), B(x_j, \eta^{- C} |I_j|^{\frac12})\bigr)
&\lesssim \bigl[ \eta^{- C} |I_j|^{\frac12}\bigr]^{\frac{4n}{n+2}} \| u_{-}(\tau_j) \|^2_{\frac{2(n+2)}{n-2}} \\
&\lesssim \eta^{2C} |I_j|
\end{align*}
provided $C_1$ is chosen large enough, for example, $C_1 =3 C$. (We also require that $C_1 < 45 C$ for use in Corollary~\ref{no soliton}.)
Using \eqref{rate mass change}  as before, we find
\begin{equation}\label{u minus antibubble}
\Mass\bigl( u_{-}(t_j), B(x_j, \eta^{-2 C} |I_j|^{\frac12})\bigr)
\lesssim \eta^{2C} |I_j|.
\end{equation}

Combining \eqref{v bubble} and \eqref{u minus antibubble} with the triangle inequality we obtain
\begin{equation}\label{u bubble}
\Mass\bigl( u(t_j), B(x_j, \eta^{- C} |I_j|^{\frac12})\bigr)  \gtrsim \eta^{C} |I_j|,
\end{equation}
which yields the claim through one further application of \eqref{rate mass change}.
\end{proof}

As in our predecessors, \cite{borg:scatter,tao: gwp radial}, the radial assumption is used to show that the bubble of
mass concentration just exhibited must occur at the spatial origin.  Technology to deal with non-radial data has recently
been developed, \cite{ckstt:gwp,RV,Monica:thesis,MV}, primarily based off the notion of a frequency localized interaction Morawetz
inequality.  As a result of the necessity of introducing frequency localization, implementing this strategy in the presence of a potential
would constitute a major undertaking.

\begin{corollary}[Bubble at the origin]\label{O bubble}
Let $I_j$ be an unexceptional interval.  Then
$$
\Mass\bigl( u(t), B(0, \eta^{-6 C} |I_j|^{\frac12} ) \bigr) \gtrsim \eta^{C} |I_j|
$$
for all $t\in I_j$.  Here we may choose $C=4 n^2$.
\end{corollary}

\begin{proof}
If $x_j$ in Proposition~\ref{bubble} is within $\tfrac12 \eta^{-6 C} |I_j|^{\frac12}$ of the origin, then the result
follows immediately.  We will now show that this must occur.  For if not, by the radial assumption there would be at least
$O(\eta^{-5(n-1)C})$ many disjoint balls each containing at least $\eta^{C} |I_j|$ amount of mass.  By H\"older's inequality, this
would imply
$$
\eta^{-5(n-1)C} \times \eta^{C} |I_j|
    \lesssim  \| u(t) \|_{\frac{2n}{n-2}}^2 \times \bigl[\eta^{-6 C} |I_j|^{\frac12} \bigr]^{\frac{2(n-1)}{n}}
            \bigl[\eta^{-C} |I_j|^{\frac12}\bigr]^{\frac2{n}},
$$
or more succinctly, $\| u(t) \|_{\frac{2n}{n-2}} \gtrsim \eta^{-\theta C}$, where $\theta = \tfrac1{2n}(5n^2-18n+10)$.
As $\theta>0$, this contradicts our \emph{a priori} bound on the potential energy (see \eqref{pot decay} in the repulsive case).
\end{proof}

\subsection{Morawetz to the rescue}
The Morawetz inequality derived in subsection~\ref{SS:morawetz} clearly speaks against the
existence of bubbles of mass concentration at the origin.  Unfortunately, it is not sufficient
by itself to limit the total number of bubbles (and so unexceptional intervals).  The way to
proceed was discovered by Bourgain, \cite{borg:scatter}; it is ingenious.

Morawetz does limit the number of nearby intervals of any given scale; this is encapsulated in Corollary~\ref{no soliton}
below. Thus, if there were many unexceptional intervals, they must form a cascade --- an accumulating sequence of
diadically shrinking intervals. This in turn is prohibited by a result that shows that energy cannot rapidly move to high frequencies;
this is usually termed `energy non-evacuation'.

\begin{corollary}\label{no soliton}
Let $u$ be a radial solution of \eqref{equation} on $I$, a compact time interval contained in $[t_-,t_+]$.  We partition $I=\cup I_j$ as
in \eqref{small W}.  Then
\begin{equation*}
\sum |I_j|^\frac12 \lesssim \eta^{-51 C} |I|^\frac12
\end{equation*}
and moreover,
\begin{equation*}
 |I_j| \gtrsim   \eta^{102 C}  |I|,
\end{equation*}
for at least one index $j$.
\end{corollary}

\begin{proof}
This result will follow directly from Corollary~\ref{O bubble} and Proposition~\ref{morawetz}.

From H\"older's inequality and Corollary~\ref{O bubble}, we have
\begin{align*}
\int_{|x|\leq R} \frac{|u(t,x)|^{\frac{2n}{n-2}}}{|x|} \, dx
&\gtrsim  R^{-1-\frac{2n}{n-2}} \bigl[ \eta^{C} |I_j| \bigr]^{\frac{n}{n-2}}
\end{align*}
and hence,
\begin{align*}
\int_{|x|\leq R} \frac{|u(t,x)|^{\frac{2n}{n-2}}}{|x|} \, dx
    \gtrsim \eta^{45C} |I_j|^{-\frac12}
\end{align*}
for any $R\geq \eta^{-6 C} |I_j|^{\frac12}$ and any unexceptional interval $I_j$.
We now integrate this over each $I_j$ and sum; combining this with Proposition~\ref{morawetz} leads to
$$
\sum |I_j|^\frac12 \lesssim \eta^{-51 C} |I|^\frac12.
$$
We do not need to restrict the sum to unexceptional intervals as the total number of exceptional intervals
is $O(\eta^{-\frac{2(n+2)}{n-2}C_1})$; see Remark~\ref{count ex}.

The second claim follows from the first by writing $|I_j|^{1/2} \geq |I_j| [ \sup |I_k| ]^{-1/2}$.
\end{proof}

The proof of the existence of a cascade of unexceptional intervals is based on the following general proposition.
This is implicit in \cite{borg:scatter} and was made explicit in \cite{tao: gwp radial}.  We present the latter proof
with a few expository modifications.

\begin{proposition}[Interval cascade]\label{cascade}
Let $I$ be an interval tiled by finitely many intervals $I_1,\ldots,I_N$.  Suppose that for any
contiguous family $\{I_j : j\in\mathcal{J}\}$ there exists $j_*\in\mathcal{J}$ so that
\begin{equation}\label{scale hyp}
|I_{j_*}| \geq a |\cup_{j\in\mathcal{J}} I_j|
\end{equation}
for some small $a>0$.  Then there exist $K\geq \log(N)/\log(2a^{-1})$ distinct indices $j_1,\ldots,j_K$ such that
$$
|I_{j_1}| \geq 2 |I_{j_2}| \geq \cdots \geq 2^{K-1} |I_{j_K}|
$$
and $I_{j_{l}} \subseteq 3 a^{-1}  I_{j_k}$ for all $l<k$.
\end{proposition}

\begin{proof}
We begin by running an algorithm that assigns a generation (a positive integer) to each $I_j$.

By hypothesis, $I$ contains at least one interval of length $a|I|$.  All intervals of length $\geq \frac12 a |I|$
belong to the first generation.  By looking at the total measure, we see that there are at most $2a^{-1}-1$ intervals
in the first generation. Removing these intervals from $I$ leaves at most $2a^{-1}$ gaps, which are of course tiled by intervals $I_j$.

Notice that by \eqref{scale hyp} no gap is longer than $\tfrac12 |I|$.  As a result, any interval $I_j$ that has not been selected
to belong to the first generation obeys $I_j\subseteq (2 a^{-1}+1) I_k$ for some $I_k$ in the first generation.

We now apply this argument recursively to the gaps generated by the previous iteration (one should view $I$ as the gaps from the
zeroth iteration) until every $I_j$ has been labeled with a generation number.

Each iteration of the algorithm removes at most $2a^{-1}-1$ many intervals and produces at most $2a^{-1}$ gaps.
As we begin with $N$ intervals, we are guaranteed that the number $K$ of iterations performed must obey
$$
N \leq  (2a^{-1}-1) + (2a^{-1}-1)2a^{-1} + \cdots + (2a^{-1}-1)[2a^{-1}]^{K-1}.
$$
A few manipulations yield the claim.
\end{proof}

\begin{proposition}[Energy non-evacuation]\label{E non-evac}
Let $I_{j_1},\ldots,I_{j_K}$ be a disjoint family of unexceptional intervals obeying
$$
|I_{j_1}| \geq 2 |I_{j_2}| \geq \cdots \geq 2^{K-1} |I_{j_K}|
$$
and $t_*$ a time such that $\dist( I_{j_{k}} , t_* ) \lesssim \eta^{-102 C}  |I_{j_k}|$ for $1\leq k \leq K$.
Then $K\leq \eta^{-10^3 C}$.
\end{proposition}

\begin{proof}
By Corollary~\ref{O bubble},
$$
\Mass\bigl( u(t), B(0, \eta^{-6 C} |I_{j_k}|^{\frac12} ) \bigr) \gtrsim \eta^{C} |I_{j_k}|
$$
for all $t\in I_{j_k}$.  Moreover, from \eqref{rate mass change} we can deduce that this mass does not move too far away
by time $t_*$.  Specifically,
$$
\Mass\bigl( u(t_*), B(0, \eta^{-103 C} |I_{j_k}|^{\frac12} ) \bigr) \gtrsim \eta^{C} |I_{j_k}|.
$$
Conversely, by \eqref{mass in ball},
$$
\Mass\bigl( u(t_*), B(0, R ) \bigr) \lesssim  R^2.
$$
Putting these together we obtain
$$
\int_{A(k)} |u(t_*,x)|^2 \, dx \gtrsim \eta^{C} |I_{j_k}|,
$$
where $A(k)$ is the annulus
$$
A(k) = \bigl\{ x :  \eta^C |I_{j_k}|^{\frac12} \leq |x| \leq \eta^{-103 C} |I_{j_k}|^{\frac12} \bigr\}.
$$

Choosing $M\simeq 105 C \log(\frac1\eta)$, the annuli associated to $k=1, M+1, 2M+1,\ldots$ are disjoint.
The number of such annuli is $O(K/M)$.

By H\"older's inequality,
\begin{align*}
\int_{A(k)} |u(t_*,x)|^\frac{2n}{n-2} \, dx
&\gtrsim \bigl[ \eta^{C} |I_{j_k}| \bigr]^\frac{n}{n-2} \bigl[ \eta^{-103 C} |I_{j_k}|^\frac12 \bigr]^{-\frac{2n}{n-2}}
    \geq \eta^{700 C}
\end{align*}
and so by using the fact disjointness of the annuli chosen above, we obtain
$$
\tfrac{K}{M} \eta^{700 C} \lesssim \|u(t_*)\|_\frac{n}{n-2}^\frac{n}{n-2} \lesssim 1.
$$
This gives the bound on $K$ we claimed --- indeed, with a lot of room to spare.
\end{proof}

We are now prepared to complete the proof of global well-posedness.  As explained in subsection~\ref{SS:6.1},
our sole obligation is to bound the total number of unexceptional intervals; by Remark~\ref{count ex}, the number of exceptional intervals
is $O(\eta^{-\frac{2(n+2)}{n-2}C_1})$.

We first bound the number $N$ of unexceptional intervals that can occur consecutively.  Let us write $I$ for the
union of these intervals.   By Corollary~\ref{no soliton}, the hypotheses of Proposition~\ref{cascade} are
satisfied with $a=\eta^{-102 C}$ and so we can find a cascade of $K$ intervals.  By choosing $t_*$ to be any
point in $I_{j_K}$ this collection of intervals will satisfy the hypotheses of Proposition~\ref{E non-evac}.
The resulting bound on $K$ induces (through Proposition~\ref{cascade}) a bound on $N$, namely
$$
N \lesssim (2\eta^{-102 C})^{\eta^{-10^3C}}.
$$

Lastly, as there are only $O(\eta^{-\frac{2(n+2)}{n-2}C_1})$ exceptional intervals (which also bounds the number
of groups of consecutive unexceptional intervals), the total number of intervals is
$$
J \lesssim \eta^{-\frac{2(n+2)}{n-2}C_1} + \eta^{-\frac{2(n+2)}{n-2}C_1} N \leq \exp\{ \eta^{-10^4 C} \}.
$$
As the $Z$-norm of $u$ on each of these $J$ intervals is $\lesssim \eta$ (cf. \eqref{small Z}), this bounds the total $Z$-norm
of $u$ on $[t_-,t_+]$.  Thus \eqref{goal} and hence global well-posedness follow.

%%%%%%%%%%%%%%%%%%%%%%%%%%%%%%%%%%%%%%%%%%%%%%%%%%%%%%%%%%%%%%%%%%%%%%%%%%%%%%%%%%%%%%%%%%%
%
%
%                                   Section
%
%
%%%%%%%%%%%%%%%%%%%%%%%%%%%%%%%%%%%%%%%%%%%%%%%%%%%%%%%%%%%%%%%%%%%%%%%%%%%%%%%%%%%%%%%%%%%

\section{Scattering in the repulsive case}\label{scattering section}

In this section, we prove scattering for radial solutions to \eqref{equation} in the repulsive case.
As we noted at the beginning of the previous section, the estimate  \eqref{goal repuls} alone is not sufficient to yield scattering.
In order to establish Strichartz estimates powerful enough to imply scattering, we have to rely on the exponential decay in time of the
potential energy, as presented in Lemma~\ref{lemma energy control}.

\begin{lemma}[Good Strichartz bounds]\label{stric bdd}
Let $V(x)=-\tfrac 12 |x|^2$ and $u_0\in \Sigma$.  Suppose that on every compact time interval $[-T,T]$,
there exists a unique solution $u$ to \eqref{equation} which obeys
\begin{align*}
\|u\|_{Z([-T,T])}\leq C(T,\|u_0\|_{\Sigma}).
\end{align*}
Then,
\begin{align*}
\| u \|_{Z(\R)} + \|H(-t)^{\frac 12} u\|_{S^0(\R)} \leq C(\|u_0\|_{\Sigma}).
\end{align*}
\end{lemma}

\begin{proof}
By time-reversal symmetry, it suffices to establish the claim for positive times only.

Let $\eta>0$ be a small constant to be chosen momentarily.  Then, there exists $T_0= T_0(\|\nabla u_0\|_2)$ such that for $T\geq T_0$,
\begin{align}\label{choose T_0}
\mathcal E_1(0) \cosh^{-2}(T)\leq \eta^{\frac{2n}{n-2}}.
\end{align}
In particular, by \eqref{pot decay},
\begin{align}\label{pot small}
\|u\|_{L_t^\infty L_x^\frac{2n}{n-2}([T_0, \infty)\times\R^n)} \lesssim \eta.
\end{align}
By Duhamel's formula, on $[T_0,\infty)$, the solution $u$ satisfies
$$
u(t)=U(t-T_0)u(T_0) - i\int_{T_0}^t U(t-s) F(u(s))\, ds.
$$
Applying Strichartz, Proposition~\ref{frac chain}, and \eqref{pot small}, we estimate
\begin{align*}
\|H(-t)^{\frac 12} u\|_{S^0([T_0, \infty))}
&\lesssim \|H(-T_0)^{\frac 12} u(T_0)\|_2 + \|H(-t)^{\frac 12} F(u)\|_{L_t^2L_x^{\frac{2n}{n+2}}} \\
&\lesssim \|H(-T_0)^{\frac 12} u(T_0)\|_2 + \|u\|_{L_t^\infty L_x^\frac{2n}{n-2}}^{\frac{4}{n-2}}\|H(-t)^{\frac 12} u\|_{L_t^2L_x^{\frac{2n}{n-2}}} \\
&\lesssim \|H(-T_0)^{\frac 12} u(T_0)\|_2 + \eta^{\frac 4{n-2}}\|H(-t)^{\frac 12} u\|_{S^0([T_0, \infty))},
\end{align*}
where all spacetime norms are on $[T_0, \infty)\times\R^n$.  Thus, taking $\eta$ sufficiently small depending only on the
Strichartz constant, we get
\begin{align*}
\|H(-t)^{\frac 12} u\|_{S^0([T_0, \infty))}\lesssim \|H(-T_0)^{\frac 12} u(T_0)\|_2.
\end{align*}
On the other hand, by hypothesis and Lemma~\ref{stric<Z},
\begin{align*}
\|H(-t)^{\frac 12} u\|_{ S^0([0,T_0])}\leq C(T_0, \|u_0\|_\Sigma).
\end{align*}
Thus,
\begin{align*}
\|H(-t)^{\frac 12} u\|_{S^0([0,\infty))}\leq C(T_0, \|u_0\|_\Sigma) \leq C(\|u_0\|_\Sigma),
\end{align*}
which yields the stated bound on $H(-t)^{\frac 12} u$.

To obtain a bound on the $Z$-norm of $u$, we first note that by replacing $H(-t)^{\frac12}$ with $P(-t)$
in the argument just presented yields
\begin{align*}
\|P(-t) u\|_{S^0([T_0, \infty))} \lesssim \|P(-T_0) u(T_0)\|_2 \lesssim \|H(-T_0)^{\frac 12} u(T_0)\|_2 \leq C(T_0, \|u_0\|_\Sigma).
\end{align*}
Applying Lemma~\ref{lemma pt sobolev}, this implies
\begin{align*}
\| u\|_{Z([T_0, \infty))} \lesssim \cosh^{-1}(T_0) \|P(-t) u\|_{S^0([T_0, \infty))} \lesssim C(T_0, \|u_0\|_\Sigma).
\end{align*}
This completes the proof since the $Z-$norm on $[0,T_0]$ is controlled by hypothesis and $T_0$ was chosen depending only
on $\|\nabla u_0\|_2$.
\end{proof}

Combining the results of Section~\ref{gwp section} with Lemma~\ref{stric bdd}, we have settled the first half of
Theorem~\ref{repulsive gwp}, that is, global existence and good Strichartz control.  Next we establish asymptotic completeness, thus
completing the proof of Theorem~\ref{repulsive gwp}.

\begin{lemma}[Asymptotic completeness]\label{asy compl}
Let $V(x)=-\tfrac 12 |x|^2$ and $u_0\in \Sigma$.  Suppose that there exists a unique global solution $u$ to \eqref{equation} such that
\begin{align}\label{asy compl ass}
\|H(-t)^{\frac 12} u\|_{ S^0(\R)} \leq C(\|u_0\|_{\Sigma}).
\end{align}
Then, there exist unique functions $u_{\pm}\in \Sigma$ such that
$$
\|U(-t)u(t)-u_{\pm}\|_{\Sigma}\to 0, \quad  \text{as } t\to \pm\infty.
$$
\end{lemma}

\begin{proof}
We will only construct $u_+$ and show that it satisfies the claim.  The proof involving $u_-$ is basically identical
and we omit it.

We start by constructing the scattering state $u_+$.  For $t>0$ define $v(t) = U(-t)u(t)$. We will show
that $v(t)$ converges in $\Sigma$ as $t\rightarrow \infty$, and define $u_{+}$ to be that limit.
Indeed, from Duhamel's formula \eqref{duhamel} we have
\begin{align}\label{v}
v(t) = u_0 - i\int_{0}^{t} U(-s)F(u(s))\, ds.
\end{align}
Therefore, for $0<\tau<t$,
$$
v(t)-v(\tau)=-i\int_{\tau}^{t}U(-s)F(u(s))\, ds.
$$
By Lemma~\ref{L:same norm}, Strichartz, Proposition~\ref{frac chain}, Lemma~\ref{sob em}, and \eqref{pot decay},
\begin{align*}
\|v(t)-v(\tau)\|_{\Sigma}
     &\simeq \|H_0^{\frac 12}[v(t)-v(\tau)]\|_2\\
     &\simeq \Bigl\| \int_{\tau}^{t}U(-s)H(-s)^{\frac 12} F(u(s)) \, ds \Bigr\|_2\\
     &\lesssim \|H(-s)^{\frac 12} F(u)\|_{L^2_tL_x^{\frac{2n}{n+2}}([\tau,t]\times\R^n)} \\
     &\lesssim \|u\|_{L_t^\infty L_x^{\frac{2n}{n-2}}([\tau,t])}^{\frac4{n-2}}\|H(-s)^{\frac 12} u\|_{S^0([\tau,t])} \\
     &\leq C(\|\nabla u_0\|_2) \|H(-s)^{\frac 12} u\|_{S^0([\tau,t])}.
\end{align*}
By \eqref{asy compl ass} and the Dominated Convergence Theorem,
$$
\|v(t)-v(\tau)\|_{\Sigma}\rightarrow 0 \quad \text{as } t,\tau\to \infty.
$$
In particular, this implies that $u_{+}$ is well defined.  Also, inspecting \eqref{v} one easily sees that
\begin{align*}
u_{+}=u_0- i\int_{0}^{\infty}U(-s)F(u(s))\, ds.
\end{align*}
Thus, arguing as before,
\begin{align*}
\|U(-t)u(t)-u_+\|_{\Sigma}
&=\Bigl\|\int_{t}^{\infty}U(-s)F(u(s))\, ds\Bigr\|_{\Sigma}\\
&\simeq \Bigl\|\int_{t}^{\infty}H_0^{\frac 12} U(-s)F(u(s))\, ds\Bigr\|_2\\
&\leq C(\|\nabla u_0\|_2) \|H(-s)^{\frac 12} u\|_{S^0([t, \infty))}\to 0 \quad \text{as } t\to \infty.
\end{align*}
This completes the proof of Lemma~\ref{asy compl}.
\end{proof}

Finally, we construct the wave operators.

\begin{lemma}[Existence of wave operators]\label{wave op}
Let $V(x)=-\tfrac 12 |x|^2$ and $u_+\in \Sigma$ be radial.  Then, there exists a unique global solution $u$ to \eqref{equation} with $u_0\in \Sigma$
such that
\begin{align}\label{wave op concl}
\|H(-t)^{\frac 12} u\|_{ S^0(\R)} \leq C(\|u_+\|_{\Sigma})
\end{align}
and moreover,
$$
\|U(-t)u(t)-u_+\|_{\Sigma}\to 0 \quad  \text{as } t\to \infty.
$$
\end{lemma}

\begin{proof}
By time reversal symmetry, a similar result holds in the negative time direction.  More precisely, for $u_-\in \Sigma$ radial
there exists a unique global solution $u$ to \eqref{equation} with $u_0\in \Sigma$ such that
\begin{align*}
\|H(-t)^{\frac 12} u\|_{ S^0(\R)} \leq C(\|u_-\|_{\Sigma})
\end{align*}
and
$$
\|U(-t)u(t)-u_-\|_{\Sigma}\to 0 \quad  \text{as } t\to -\infty.
$$
As the proof is basically identical, we only present it for the positive time direction.

By standard arguments, it suffices to show that the integral equation
\begin{align*}
u(t)=U(t)u_+ + i\int_{t}^{\infty}U(t-s)F(u(s))\, ds
\end{align*}
admits a unique global solution that satisfies \eqref{wave op concl}.

First, we observe that for $T=T(\|\nabla u_+\|_2)$ sufficiently large, the solution map
\begin{align*}
[\Phi(u)](t)=U(t)u_+ + i\int_{t}^{\infty}U(t-s)F(u(s))\, ds
\end{align*}
is a contraction on the ball
\begin{align*}
\mathcal B=\{u;\ & \|P(-t) u\|_{L_t^{\frac{2(n+2)}{n-2}}L_x^{\frac{2n(n+2)}{n^2+4}}([T,\infty)\times\R^n)}\le 2 C \|\nabla u_+\|_2, \\
                 & \|X(-t) u\|_{L_t^{\frac{2(n+2)}{n-2}}L_x^{\frac{2n(n+2)}{n^2+4}}([T,\infty)\times\R^n)}\le 2 C \|xu_+\|_2,\\
                 & \|u\|_{L_t^{\frac{2(n+2)}{n-2}}L_x^{\frac{2n(n+2)}{n^2+4}}([T,\infty)\times\R^n)}\le 2 C \|u_+\|_2\}
\end{align*}
equipped with the distance given by the $L_t^{\frac{2(n+2)}{n-2}}L_x^{\frac{2n(n+2)}{n^2+4}}$-norm.
Here, $C$ is the Strichartz constant.  The proof is immediate and follows from Strichartz estimates; see Proposition~\ref{lwp vn0} for more details.
We only mention that the smallness necessary to close the argument comes from the Sobolev embedding inequality
(see Lemma~\ref{lemma pt sobolev})
$$
\|u\|_{Z([T,\infty))} \lesssim \cosh^{-1}(T)\|P(-t)u\|_{L_t^{\frac{2(n+2)}{n-2}}L_x^{\frac{2n(n+2)}{n^2+4}}([T,\infty)\times\R^n)}
$$
and taking $T$ sufficiently large (depending only on $\|\nabla u_+\|_2$).

Note that $u(T)\in \Sigma$ and morover,
$$
\|u(T)\|_\Sigma \leq C(\|u_+\|_\Sigma).
$$
At this point, using the results of Section~\ref{gwp section} we can solve the Cauchy problem with initial data $u(T)$ at time $T$.
The unique solution satisfies
$$
\|H(-t)^{\frac 12} u\|_{S^0((-\infty,T])} \leq C(T, \|u(T)\|_{\Sigma}) \leq C(\|u_+\|_\Sigma).
$$
Putting everything together, we derive \eqref{wave op concl}.  This concludes the proof of Lemma~\ref{wave op}.
\end{proof}

\end{document}